\documentclass[a4paper, 10pt, twoside, notitlepage]{amsart}

\usepackage{amsmath,amscd}
\usepackage{amssymb}
\usepackage{amsthm}
\usepackage{comment}
\usepackage{graphicx, xcolor}

\usepackage{mathrsfs}
\usepackage[ocgcolorlinks, linkcolor=blue]{hyperref}

\usepackage{bm}
\usepackage{bbm}
\usepackage{url}

\usepackage[utf8]{inputenc}
\usepackage{mathtools,amssymb}
\usepackage{esint}
\usepackage{tikz}
\usepackage{dsfont}
\usepackage{relsize}
\usepackage{url}
\usepackage{xcolor}
\usepackage{graphicx}
\usepackage{mathrsfs}
\usepackage[shortlabels]{enumitem}
\usepackage{lineno}
\usepackage{amsmath}
\usepackage{enumitem}
\usepackage{amsthm} 
\usepackage{verbatim}
\usepackage{dsfont}
\numberwithin{equation}{section}

\allowdisplaybreaks

\mathtoolsset{showonlyrefs}

\graphicspath{{images/}}

\newtheorem{theorem}{Theorem}[section]
\newtheorem{lemma}[theorem]{Lemma}
\newtheorem{definition}[theorem]{Definition}
\newtheorem{corollary}[theorem]{Corollary}

\newtheorem{proposition}[theorem]{Proposition}
\newtheorem{example}[theorem]{Example}

\newtheorem{remark}[theorem]{Remark}

{\theoremstyle{remark}
\newtheorem{claim}{Claim}}

\author[P. Zimmermann]{Philipp Zimmermann}
\address{Institute of Mathematics, École Polytechnique Fédérale de Lausanne (EPFL), Lausanne, Switzerland}
\email{philipp.zimmermann@epfl.ch}

\subjclass[2020]{35R30, 26A33, 42B37}
\keywords{Viscous wave equations, nonlocal operators, inverse problems, unique continuation principle, Runge approximation theorem}

\newcommand{\C}{{\mathbb C}}
\newcommand{\R}{{\mathbb R}}
\newcommand{\Z}{{\mathbb Z}}
\newcommand{\N}{{\mathbb N}}

\newcommand{\eps}{\varepsilon}

\newcommand{\re}{\text{ Re}} 

\newcommand{\schwartz}{\mathscr{S}}

\newcommand{\test}{\mathscr{D}}

\DeclareMathOperator{\supp}{supp} 
\DeclareMathOperator{\dist}{dist} 

\newcommand{\weak}{\rightharpoonup}
\newcommand{\weakstar}{\rightharpoonup^\ast}

\title[Weakly viscous nonlocal waves]{Weakly viscous nonlocal wave equations: forward theory, Lebesgue Runge approximation, and inverse problems}
\begin{document}
\maketitle
\begin{abstract}
We study forward and inverse problems for weakly viscous nonlocal wave equations with elastic term $(-\Delta)^s u$ and dissipative term $(-\Delta)^r\partial_tu$, where $0\leq r<s$. For time-dependent potentials in suitable $L^\infty(0,T;L^p(\Omega))$ classes, we establish finite energy well-posedness, strong continuity in time, and an energy identity by parabolic regularization and time-regularization arguments. Under explicit nonlocality and nonresonance assumptions on the fractional orders, we prove Runge approximation in $L^2(0,T;\widetilde H^r(\Omega))$ and hence in $L^2(\Omega\times(0,T))$. As an application, partial exterior Dirichlet-to-Neumann measurements uniquely determine the potential in all three exponent regimes of the Runge theorem, including integer elastic orders with noninteger dissipative orders. The recovery argument combines equality of solutions with the same exterior input, unique continuation or the entanglement principle, and a single Lebesgue-space Runge approximation. It applies to the full admissible potential class, without spatial boundedness, sign, or smallness assumptions.
\end{abstract}

	\tableofcontents

\section{Introduction}
\label{intro:section}

The recovery of an interior coefficient from exterior measurements is a
central question in inverse problems for nonlocal equations. For fractional
wave equations, the interaction between nonlocality and time evolution makes
it possible to probe the interior by prescribing data on an arbitrarily
chosen open subset of the exterior. When dissipation is also nonlocal, the
measurement contains contributions from both the displacement and its time
derivative. Understanding how to separate these contributions requires a
forward theory that respects their different spatial regularities.

In this article, we study the forward problem, Runge approximation, and
potential recovery for weakly viscous nonlocal wave equations of the form
\begin{equation}\label{intro:equation}
\begin{cases}
\bigl(\partial_t^2+(-\Delta)^r\partial_t+(-\Delta)^s+q\bigr)u=f
    &\text{in }\Omega_T,\\
u=\varphi&\text{in }(\Omega_e)_T,\\
u(0)=u_0,\quad\partial_tu(0)=u_1&\text{in }\Omega,
\end{cases}
\end{equation}
where $\Omega\subset\R^n$ is a bounded Lipschitz domain,
$\Omega_e=\R^n\setminus\overline\Omega$, $T>0$, and
\[
0\leq r<s<\infty.
\]
For a set $A\subset\R^n$, we write $A_T=A\times(0,T)$.
The fractional Laplacian $(-\Delta)^\rho$ is the Fourier multiplier with
symbol $|\xi|^{2\rho}$, and $q=q(x,t)$ is a real-valued potential.
We first construct finite energy solutions for the full range $0\leq r<s$.
Under additional nonlocality and nonresonance assumptions, we then prove
Runge approximation and use it to determine $q$ from partial exterior data.
These three results are connected by the same energy framework: the forward
theory supplies the adjoint solutions needed for approximation, and
approximation supplies the interior test functions needed for recovery.

\subsection{Weakly viscous waves and related models}

The terminology \emph{weakly viscous} refers here to the inequality $r<s$:
the spatial order of the dissipation is smaller than that of the elastic
term. The endpoint $r=0$ gives the damped equation
\[
\partial_t^2u+\partial_tu+(-\Delta)^su+qu=f,
\]
which the author investigated in \cite{zimmermann2024optimalrungeapproximationdamped}, whereas $r=s$, outside the range considered here, gives the nonlocal viscous
wave equation studied in \cite{zimmermann2024-viscous-wave}.
Allowing distinct powers separates the length scales on which the elastic
response and dissipation act. A concrete example is
\[
\partial_t^2u+(-\Delta)^{1/2}\partial_tu-\Delta u=0,
\]
whose nonlinear counterparts arise in models of fluid--structure interaction;
see \cite{kuan2021deterministic,kuan2022probabilistic}.
Here the elastic term is local, but dissipation is nonlocal. Such examples
also motivate treating integer and noninteger elastic orders within the
same forward theory.

For homogeneous exterior data, vanishing forcing and $q=0$, the formal
energy balance is
\begin{equation}\label{intro:energy}
\frac{d}{dt}\left(
\|\partial_tu(t)\|_{L^2(\Omega)}^2+
\|(-\Delta)^{s/2}u(t)\|_{L^2(\R^n)}^2\right)
=-2\|(-\Delta)^{r/2}\partial_tu(t)\|_{L^2(\R^n)}^2.
\end{equation}
Thus dissipation controls the velocity in $H^r$, while the elastic energy
controls the displacement in $H^s$. When $r<s$, these controls do not give
$\partial_tu\in L^2(0,T;H^s)$. In particular, testing the elastic term
directly against $\partial_tu$ requires a justification beyond the
regularity initially available for weak solutions.

We address this difficulty by adding
$\varepsilon(-\Delta)^s\partial_tu_\varepsilon$ to the equation.
The regularized problems fall within an abstract viscous evolution
framework. Uniform estimates yield a weak limit as $\varepsilon\to0$;
a further time-regularization argument establishes the energy identity
and strong continuity in the energy spaces. An integrated-in-time test
function gives uniqueness in the larger class of weak finite energy
solutions. This construction accommodates time-dependent potentials
without imposing a sign or smallness condition.

\subsection{Main results}

We first specify the potential class. As in
Definition~\ref{F:class of potentials}, $\mathcal Q_s$ consists of the
real-valued potentials that belong to $L^\infty(0,T;L^p(\Omega))$ for some
$2\leq p\leq\infty$ satisfying
\begin{equation}\label{intro:potential-exponents}
\begin{cases}
n/s\leq p\leq\infty,&2s<n,\\
2<p\leq\infty,&2s=n,\\
2\leq p\leq\infty,&2s>n,
\end{cases}
\end{equation}
and are distributionally continuous in time, meaning that
$t\mapsto\int_\Omega q(x,t)\zeta(x)\,dx$ is continuous on $[0,T]$
for every $\zeta\in C_c^\infty(\Omega)$.
We use $\widetilde H^\rho(\Omega)$ for the closure of
$C_c^\infty(\Omega)$ in $H^\rho(\R^n)$, and $H^{-\rho}(\Omega)$ for
its antidual. The space $\widetilde L^2(\Omega)$ consists of all $L^2(\R^n)$ functions that vanish outside $\Omega$ and it can be identified with
$L^2(\Omega)$.

\begin{theorem}[Finite energy well-posedness]\label{intro:forward}
Let $0\leq r<s$, $q\in\mathcal Q_s$, and
\[
(u_0,u_1,f)\in\widetilde H^s(\Omega)\times\widetilde L^2(\Omega)
\times L^2(0,T;H^{-r}(\Omega)).
\]
Then \eqref{intro:equation} with $\varphi=0$ has a unique weak finite
energy solution. This solution belongs to
\[
C([0,T];\widetilde H^s(\Omega))
\cap C^1([0,T];\widetilde L^2(\Omega))
\cap H^1(0,T;\widetilde H^r(\Omega)),
\]
attains the prescribed initial data, and satisfies the energy identity
\eqref{F:eq: energy identity weakly viscous nonlocal wave eq continuous potential}.
\end{theorem}

This statement combines Theorems~\ref{F:prop: existence of weak solutions}
and~\ref{F:prop: uniqueness}. Nonzero exterior data are handled by lifting;
see Theorem~\ref{F:thm: unique solvability of weakly viscous wave eq nonhom}.
In particular, it applies to the smooth exterior inputs with zero initial
data used in the inverse problem below. No noninteger-order assumption is
needed for the forward result.

Our approximation theorem does require nonlocality. Given a nonempty open
set $W\subset\Omega_e$, let $u_\varphi$ solve \eqref{intro:equation} with
$f=u_0=u_1=0$ and exterior datum $\varphi\in C_c^\infty(W_T)$, and set
\[
\mathscr R_W=\{u_\varphi-\varphi:\varphi\in C_c^\infty(W_T)\}.
\]
We use the following three regimes:
\begin{enumerate}[(i)]
\item\label{intro:regime-i} $s\in\N$ and $r\notin\N_0$;
\item\label{intro:regime-ii} $s\notin\N$ and $r\in\N_0$;
\item\label{intro:regime-iii} $n\geq2$, $r,s\notin\N_0$, and
\begin{equation}\label{intro:nonresonance}
\begin{cases}
s-r\notin\N,&n\text{ even},\\
s-r\notin\tfrac12\N,&n\text{ odd}.
\end{cases}
\end{equation}
\end{enumerate}
Throughout these alternatives, $0\leq r<s$ remains understood.

\begin{theorem}[Runge approximation]\label{intro:runge}
Let $q\in\mathcal Q_s$ and suppose that one of the regimes
\ref{intro:regime-i}--\ref{intro:regime-iii} holds.
For every nonempty open set $W\subset\Omega_e$, the Runge set
$\mathscr R_W$ is dense in $L^2(0,T;\widetilde H^r(\Omega))$.
Consequently, the restrictions $u_\varphi|_{\Omega_T}$ are dense in
$L^2(\Omega_T)$.
\end{theorem}

Theorem~\ref{R:thm: Runge approx in Hr} proves the spatial Sobolev
approximation, and Corollary~\ref{new:L2} gives its Lebesgue-space
consequence. The proof uses Hahn--Banach duality and a terminal-value
adjoint equation. The forward theory permits sources in
$L^2(0,T;H^{-r}(\Omega))$, precisely the dual space needed for the
stated approximation topology. Orthogonality to the Runge set produces
an exterior identity involving two fractional powers. When one power
is integer, locality reduces the identity to a unique continuation
problem. When both are noninteger, the entanglement principle separates
them under \eqref{intro:nonresonance}.

To state the inverse result, define the \emph{Dirichlet-to-Neumann map} by
\begin{equation}
\label{intro:DN}
\Lambda_q\varphi
=\left.\bigl((-\Delta)^r\partial_t+(-\Delta)^s\bigr)
 u_{q,\varphi}\right|_{(\Omega_e)_T},
\end{equation}
where $u_{q,\varphi}$ is the solution with zero source and initial data.
The restriction is distributional; the corresponding energy pairing
is given in \eqref{new:DN-pairing}.

\begin{theorem}[Partial-data recovery of a time-dependent potential]
\label{intro:inverse}
Suppose that one of the regimes~\ref{intro:regime-i}--\ref{intro:regime-iii} holds.
Let $q_1,q_2\in\mathcal Q_s$, and let $W_1,W_2\subset\Omega_e$ be
nonempty open sets. If
\[
\left.\Lambda_{q_1}\varphi\right|_{(W_2)_T}
=\left.\Lambda_{q_2}\varphi\right|_{(W_2)_T}
\quad\text{for every }\varphi\in C_c^\infty((W_1)_T),
\]
then $q_1=q_2$ almost everywhere in $\Omega_T$.
\end{theorem}

This is Theorem~\ref{new:potential} together with
Remark~\ref{new:unbounded-potential}. The measurement sets may overlap,
and the potentials need not be bounded, nonnegative, small, or invariant
under time reversal. The inverse theorem covers the same exponent regimes
as the Runge theorem. When $s$ is an integer and $r$ is noninteger,
unique continuation first determines the difference of the velocities;
the common initial displacement then determines the difference of the solutions.

The recovery argument compares the two solutions for the same exterior
input. Equality of the measurements and unique continuation or
entanglement give equality of these solutions throughout the domain.
Subtracting their equations therefore yields
$(q_1-q_2)u_{q_1,\varphi}=0$. Approximate the constant function $1$ in
$L^2(\Omega_T)$ by such solutions. Since
$q_1-q_2\in L^2(\Omega_T)$, a single Cauchy--Schwarz estimate passes the
identity to the limit in $L^1$. Thus Lebesgue approximation suffices even
for the unbounded potentials in $\mathcal Q_s$. This step explains why
the inverse problem can be completed within the finite energy theory
of the present article.

\subsection{Relation to inverse problems for nonlocal waves}

The fractional Calder\'on problem studied by Ghosh, Salo, and Uhlmann
\cite{GSU20} established the importance of unique continuation and
Runge approximation for recovery from exterior measurements.
For the fractional wave equation, Kow, Lin, and Wang
\cite{KLW2022} developed partial-data uniqueness and stability results.
These works provide the basic connection between exterior observations
and interior approximation that also underlies our analysis.

For the viscous equation with coincident powers $r=s$,
\cite{zimmermann2024-viscous-wave} proves recovery of linear potentials
and homogeneous nonlinearities, together with spatial-energy Runge
approximation. In its proof for general time-dependent potentials,
equality of the measurements first gives equality of the solutions with
the same input. We retain that recovery strategy, but the distinct
powers in \eqref{intro:equation} require a different exterior argument.
The multi-term entanglement principle of Feizmohammadi and Lin
\cite{feizmohammadi2024entanglementprinciplefractionallaplacian}
provides this argument in the nonresonant fractional regimes.

A related direction pursued by Kian and Uhlmann is
\cite{kian2026simultaneousrecoveryvarioustimedependent}, who study
hyperbolic equations with space--time nonlocal attenuation. 
They establish simultaneous recovery of classes of leading and
lower-order coefficients depending on both space and time from a
local source-to-solution map, including configurations in which
the source and observation sets are disjoint in both space and time.
Their results illustrate how nonlocal attenuation can provide
additional information for inverse coefficient problems.
The present work considers spatially nonlocal dissipation of the
form $(-\Delta)^r\partial_tu$ and exterior Dirichlet-to-Neumann
measurements, linking finite energy well-posedness and Runge
approximation to the recovery of a time-dependent potential. The present work and
\cite{kian2026simultaneousrecoveryvarioustimedependent}
were developed independently and contemporaneously. 

Stronger spatial approximation for nonlocal wave equations and its
applications to nonlinear inverse problems are developed in
\cite{Optimal-Runge-nonlocal-wave}; corresponding damped problems are
treated in \cite{zimmermann2024optimalrungeapproximationdamped}.
Those results use very weak solutions to accommodate rough adjoint
sources. Here the lower-order fractional dissipation already supplies
the velocity estimate needed for sources in $H^{-r}(\Omega)$, which
leads to approximation in $L^2(0,T;\widetilde H^r(\Omega))$.
The contribution of this paper is to connect this forward theory for
unequal elastic and dissipative orders with a Runge theorem and a
partial-data potential recovery result under explicit exponent
conditions. Approximation involving stronger control of the time
derivative is pursued separately in \cite{zimmermann-paperII}.

\subsection{Organization of the article}

Section~\ref{F:sec: preliminaries} introduces the Sobolev, Bochner, and
energy spaces. Section~\ref{F:sec: weak finite energy sol} develops
parabolic regularization, proves the finite energy well-posedness
results, and treats nonhomogeneous exterior data.
Section~\ref{R:sec: UCP and entanglement} records the unique continuation
and entanglement tools. Section~\ref{R:sec: Runge for constant coeff}
establishes spatial Sobolev Runge approximation, and
Section~\ref{sec:Lebesgue-Runge-inverse} derives its Lebesgue consequence
and proves the inverse result. The appendices contain the auxiliary
functional analytic arguments, the proof of the energy identity, and
the extension of the entanglement principle to the stated Bessel
potential setting.

\section{Functional analytic preliminaries}
\label{F:sec: preliminaries}

In this section, we collect the functional analytic tools used throughout the article. In particular, we introduce fractional Sobolev spaces, recall basic properties of the fractional Laplacian, and discuss vector-valued function spaces adapted to the study of weakly viscous nonlocal wave equations.

Throughout the article, all Hilbert and Banach spaces are assumed to be complex. Moreover, scalar products are taken to be linear in the first argument and antilinear in the second. If $X$ is a complex Banach space, then $X'$ denotes its antidual space, that is, the space of all continuous antilinear functionals on $X$.

\subsection{Fractional Sobolev spaces}
\label{F:sec: Fractional Sobolev spaces}

We denote by $\schwartz(\R^n)$ the Schwartz space and by $\schwartz'(\R^n)$ the space of tempered distributions. For $u\in \schwartz(\R^n)$, the Fourier transform is defined by
\[
\widehat{u}(\xi)
\vcentcolon=
\mathcal{F}u(\xi)
=
\int_{\R^n}u(x)e^{-ix\cdot \xi}\,dx.
\]
Its inverse is denoted by $\mathcal{F}^{-1}$.

For $s\in\R$ and $1<p<\infty$, the \emph{Bessel potential space} $H^{s,p}(\R^n)$ is defined by
\[
H^{s,p}(\R^n)
\vcentcolon=
\left\{
u\in \schwartz'(\R^n)
\,;\,
\langle D\rangle^s u\in L^p(\R^n)
\right\},
\]
where $\langle D\rangle^s$ denotes the Fourier multiplier with symbol
\[
\langle \xi\rangle^s
=
(1+|\xi|^2)^{s/2}.
\]
The corresponding norm is given by
\[
\|u\|_{H^{s,p}(\R^n)}
=
\|\langle D\rangle^s u\|_{L^p(\R^n)}.
\]
For $s\geq 0$, the fractional Laplacian $(-\Delta)^s$ is defined as the Fourier multiplier
\begin{equation}
\label{F:eq:fracLapFourDef}
(-\Delta)^s u
=
\mathcal{F}^{-1}
\left(
|\xi|^{2s}\widehat{u}(\xi)
\right).
\end{equation}
By the Mikhlin multiplier theorem,
\[
(-\Delta)^s
\colon
H^{r,p}(\R^n)
\to
H^{r-2s,p}(\R^n)
\]
defines a bounded linear operator for all $r\in\R$ and $1<p<\infty$.

For open sets $\Omega\subset\R^n$ and closed sets $F\subset\R^n$, we define
\begin{align*}
H^{s,p}(\Omega)
&\vcentcolon=
\left\{
u|_\Omega
\,;\,
u\in H^{s,p}(\R^n)
\right\},
\\
\widetilde{H}^{s,p}(\Omega)
&\vcentcolon=
\overline{C_c^\infty(\Omega)}^{\,H^{s,p}(\R^n)},
\\
H_F^{s,p}
&\vcentcolon=
\left\{
u\in H^{s,p}(\R^n)
\,;\,
\supp(u)\subset F
\right\}.
\end{align*}
The space $H^{s,p}(\Omega)$ is equipped with the quotient norm
\[
\|u\|_{H^{s,p}(\Omega)}
=
\inf
\left\{
\|U\|_{H^{s,p}(\R^n)}
\,;\,
U|_\Omega=u
\right\}.
\]

Throughout this article, we use the convention $H^s(\R^n)
=
H^{s,2}(\R^n)$, and similarly for the local spaces.

We also require the following fractional Poincar\'e inequality.

\begin{proposition}[{Poincar\'e inequality,\cite[Lemma~5.4]{RZ-unbounded}}]
\label{F:prop:Poincare ineq}
Let $\Omega\subset\R^n$ be a bounded domain, $0\leq r\leq s$, and $1<p<\infty$. Then there exists a constant $C=C(n,p,s,r,\Omega)>0$ such that
\begin{equation}
\label{F:eq: poincare ineq}
\|(-\Delta)^{r/2}u\|_{L^p(\R^n)}
\leq
C
\|(-\Delta)^{s/2}u\|_{L^p(\R^n)}
\end{equation}
for all $u\in C_c^\infty(\Omega)$.
\end{proposition}

By the norm equivalence
\[
\|u\|_{H^{s,p}(\R^n)}
\sim
\|u\|_{L^p(\R^n)}
+
\|(-\Delta)^{s/2}u\|_{L^p(\R^n)},
\]
Proposition~\ref{F:prop:Poincare ineq}, and the density of $C_c^\infty(\Omega)$ in $\widetilde{H}^{s,p}(\Omega)$, we obtain the following characterization.

\begin{lemma}
\label{F:lemma: equivalent norm on tilde spaces}
Let $\Omega\subset\R^n$ be a bounded domain, $s\geq 0$, and $1<p<\infty$. Then
\[
\|u\|_{\widetilde{H}^{s,p}(\Omega)}
\vcentcolon=
\|(-\Delta)^{s/2}u\|_{L^p(\R^n)}
\]
defines a norm equivalent to the standard $H^{s,p}(\R^n)$-norm on $\widetilde{H}^{s,p}(\Omega)$.
\end{lemma}

\subsection{Vector-valued spaces}
\label{F:subsec: vector-valued spaces}

We next recall several vector-valued function spaces used throughout the paper. Let $X$ be a Banach space and let $-\infty\leq a<b\leq \infty$.

For $1\leq p\leq \infty$, we denote by $L^p(a,b;X)$ the space of strongly measurable functions $u \colon (a,b)\to X$ such that $t\mapsto \|u(t)\|_X$ belongs to $L^p(a,b)$. The norm is given by
\[
\|u\|_{L^p(a,b;X)}
=
\left(
\int_a^b
\|u(t)\|_X^p\,dt
\right)^{1/p},
\]
with the usual modification for $p=\infty$.

If $a,b\in\R$ and $k\in\N_0$, we write $C^k([a,b];X)$ for the space of $k$-times continuously differentiable functions endowed with the norm
\[
\|u\|_{C^k([a,b];X)}
=
\max_{0\leq \ell\leq k}
\|\partial_t^\ell u\|_{L^\infty(a,b;X)}.
\]

Whenever no confusion can arise, we identify functions $u\in L^q(a,b;L^p(\Omega))$ with measurable functions on $\Omega\times(a,b)$ and write $u=u(x,t)$. In particular,
\[
L^p(0,T;L^p(\Omega))
=
L^p(\Omega_T).
\]

We shall also make use of vector-valued Sobolev spaces. If $X$ is a Hilbert space and $k\in\N$, we define
\begin{equation}
\label{F:eq: vector hilbert spaces}
H^k(a,b;X)
\vcentcolon=
\left\{
u\in L^2(a,b;X)
\,;\,
\partial_t^\ell u\in L^2(a,b;X)
\text{ for all }1\leq \ell\leq k
\right\}.
\end{equation}
This space becomes a Hilbert space with respect to the inner product
\[
\langle u,v\rangle_{H^k(a,b;X)}
=
\sum_{\ell=0}^k
\langle
\partial_t^\ell u,
\partial_t^\ell v
\rangle_{L^2(a,b;X)}.
\]

Finally, if $u\in L^1_{\mathrm{loc}}(a,b;X)$,
then its time derivative is understood in the distributional sense and is denoted interchangeably by
\[
\partial_t u,
\quad
u',
\quad
\frac{du}{dt}.
\]

\subsection{Energy spaces}
\label{F:subsec: energy spaces}

We now introduce the space--time energy spaces that will be used to formulate weak solutions. Throughout this subsection, let $\Omega\subset\R^n$ be a bounded Lipschitz domain, let $T>0$, and let $0\leq r\leq s$.

\begin{definition}[Energy spaces]
\label{F:def: energy spaces}
We define
\begin{equation}
\label{F:eq: intersection space}
X^1_{s,r}(\Omega_T)
\vcentcolon=
L^2(0,T;\widetilde{H}^s(\Omega))
\cap
H^1(0,T;\widetilde{H}^r(\Omega)),
\end{equation}
endowed with the Hilbert norm induced by the inner product
\begin{equation}
\label{F:eq: intersection space norm}
\langle u,v\rangle_{X^1_{s,r}(\Omega_T)}
\vcentcolon=
\langle u,v\rangle_{L^2(0,T;\widetilde{H}^s(\Omega))}
+
\langle \partial_t u,\partial_tv\rangle_{L^2(0,T;\widetilde{H}^r(\Omega))}.
\end{equation}
We also introduce the terminal-value subspace
\begin{equation}
\label{F:eq: subspace}
X^1_{s,r,T}(\Omega_T)
\vcentcolon=
L^2(0,T;\widetilde{H}^s(\Omega))
\cap
H^1_T(0,T;\widetilde{H}^r(\Omega)),
\end{equation}
where
\begin{equation}
\label{F:eq: vanishing at T}
H^1_T(0,T;\widetilde{H}^r(\Omega))
\vcentcolon=
\{u\in H^1(0,T;\widetilde{H}^r(\Omega))\,;\,u(T)=0\}.
\end{equation}
\end{definition}

The space $X^1_{s,r}(\Omega_T)$ is the natural energy space for weakly viscous equations in the regime $r < s$. Indeed, finite energy solutions typically satisfy
\[
u\in L^2(0,T;\widetilde{H}^s(\Omega)),
\quad
\partial_t u\in L^2(0,T;\widetilde{H}^r(\Omega)).
\]
The subspace $X^1_{s,r,T}(\Omega_T)$ will be used for adjoint problems and transposition formulations in the forthcoming works \cite{zimmermann-paperII,zimmermann-paperIII}, where terminal conditions are imposed at the final time horizon $T$.

\begin{remark}
\label{F:rem: equivalent description of intersection space}
In the notation of Appendix~\ref{F:appendix: generlaized vector-valued Sobolev}, one has
\[
X^1_{s,r}(\Omega_T)
=
W^1(0,T;\widetilde{H}^s(\Omega),\widetilde{H}^r(\Omega)).
\]
More generally, for $m\in\N$ one could consider
\[
X^m_{s,r}(\Omega_T)
\vcentcolon=
L^2(0,T;\widetilde{H}^s(\Omega))
\cap
H^m(0,T;\widetilde{H}^r(\Omega)).
\]
By the intermediate derivatives theorem, this coincides with the corresponding generalized vector-valued Sobolev space; see \cite[Chapter~XVIII]{DautrayLionsVol5} and Appendix~\ref{F:appendix: generlaized vector-valued Sobolev}.
\end{remark}

The next lemma collects the basic structural properties of these spaces.

\begin{lemma}[Properties of the energy spaces]
\label{F:lemma: space of sources}
The spaces $X^1_{s,r}(\Omega_T)$ and $X^1_{s,r,T}(\Omega_T)$ have the following properties.

\begin{enumerate}[(i)]
\item
\label{F:Banach prop of X}
Both $X^1_{s,r}(\Omega_T)$ and $X^1_{s,r,T}(\Omega_T)$ are Hilbert spaces. Moreover, the norm of $X^1_{s,r}(\Omega_T)$ is equivalent to the standard intersection norm
\[
\max\big(\|u\|_{L^2(0,T;\widetilde{H}^s(\Omega))},\|u\|_{H^1(0,T;\widetilde{H}^r(\Omega))}\big).
\]

\item
\label{F:Dual of X}
The antidual spaces can be identified as
\begin{equation}
\label{F:eq: dual space}
(X^1_{s,r}(\Omega_T))'
=
L^2(0,T;H^{-s}(\Omega))
+
\bigl(H^1(0,T;\widetilde{H}^r(\Omega))\bigr)'
\end{equation}
and
\begin{equation}
\label{F:eq: dual space X T}
(X^1_{s,r,T}(\Omega_T))'
=
L^2(0,T;H^{-s}(\Omega))
+
\bigl(H^1_T(0,T;\widetilde{H}^r(\Omega))\bigr)'.
\end{equation}
The sums are understood as sums of compatible Banach spaces and are endowed with the norms
\begin{equation}
\label{F:eq: norm dual X}
\inf_{F=F_1+F_2}
\left(
\|F_1\|_{L^2(0,T;H^{-s}(\Omega))}
+
\|F_2\|_{\bigl(H^1(0,T;\widetilde{H}^r(\Omega))\bigr)'}
\right)\sim \|F\|_{(X^1_{s,r}(\Omega_T))'}
\end{equation}
and
\begin{equation}
\label{F:eq: norm dual X T}
\inf_{F=F_1+F_2}
\left(
\|F_1\|_{L^2(0,T;H^{-s}(\Omega))}
+
\|F_2\|_{\bigl(H^1_T(0,T;\widetilde{H}^r(\Omega))\bigr)'}
\right)\sim \|F\|_{(X^1_{s,r,T}(\Omega_T))'}.
\end{equation}

\item
\label{F:Density in dual of X}
The spaces
\[
X^1_{s,r}(\Omega_T),
\qquad
L^2(0,T;\widetilde{L}^2(\Omega)),
\qquad
(X^1_{s,r}(\Omega_T))'
\]
form a Gelfand triple, where the second embedding is induced by the $L^2$ inner product and the identification of a Hilbert space with its antidual via the Riesz representation theorem. In other words, the continuous embeddings
\begin{equation}
\label{F:eq: dense embedding}
X^1_{s,r}(\Omega_T)
\hookrightarrow
L^2(0,T;\widetilde{L}^2(\Omega))
\hookrightarrow
(X^1_{s,r}(\Omega_T))'
\end{equation}
are dense.

\item
\label{F:density in X T}
Let $(\varphi_k)_{k\in\N}\subset \widetilde{H}^s(\Omega)$ be the eigenfunctions of the fractional Laplacian $(-\Delta)^s$ with homogeneous exterior condition. Then the set
\begin{equation}
\label{F:eq: dense set D T}
\mathcal{D}_T
\vcentcolon=
\left\{
\sum_{k=1}^N \eta_k(t)\varphi_k(x)
\,;\,
\eta_k\in C_c^\infty([0,T)),\
1\leq k\leq N,\
N\in\N
\right\}
\end{equation}
is dense in $X^1_{s,r,T}(\Omega_T)$.
\end{enumerate}
\end{lemma}

\begin{remark}
The assertions of Lemma~\ref{F:lemma: space of sources} remain valid in the abstract setting
\[
W^1_T(0,T;V,H)
\vcentcolon=
L^2(0,T;V)\cap H^1_T(0,T;H),
\]
where $V$ and $H$ are separable Hilbert spaces with dense and continuous embedding $V\hookrightarrow H$. Moreover, the density statement in Lemma~\ref{F:lemma: space of sources}\,\ref{F:density in X T} remains true if the eigenfunctions $(\varphi_k)_{k\in\N}$ are replaced by any countable dense subset of $\widetilde{H}^s(\Omega)$.
\end{remark}

\begin{proof}
The space $X^1_{s,r}(\Omega_T)$ is the intersection of two Hilbert spaces continuously embedded into the common ambient space $L^2(0,T;\widetilde{H}^r(\Omega))$. Hence it is a Hilbert space with respect to the inner product \eqref{F:eq: intersection space norm}. Since $\widetilde{H}^s(\Omega)\hookrightarrow \widetilde{H}^r(\Omega)$ continuously, the induced norm is equivalent to the standard intersection norm. This proves the assertion for $X^1_{s,r}(\Omega_T)$ in \ref{F:Banach prop of X}. The corresponding assertion for $X^1_{s,r,T}(\Omega_T)$ follows from the trace theorem, since $H^1_T(0,T;\widetilde{H}^r(\Omega))$ is a closed subspace of $H^1(0,T;\widetilde{H}^r(\Omega))$.

We next prove the antidual space identifications. Since
\[
X^1_{s,r}(\Omega_T)
=
L^2(0,T;\widetilde{H}^s(\Omega))
\cap
H^1(0,T;\widetilde{H}^r(\Omega)),
\]
with dense continuous embeddings into each factor, the standard antiduality formula for intersections of compatible complex Banach spaces gives
\[
(X^1_{s,r}(\Omega_T))'
=
L^2(0,T;H^{-s}(\Omega))
+
\bigl(H^1(0,T;\widetilde{H}^r(\Omega))\bigr)'.
\]
The norm equivalence \eqref{F:eq: norm dual X} follows from  \cite[Theorem~2.7.1]{Interpolation-spaces-Bergh-Lofstrom}. The proof of \eqref{F:eq: dual space X T} and \eqref{F:eq: norm dual X T} is identical, using the closed subspace $H^1_T(0,T;\widetilde{H}^r(\Omega))$ in place of $H^1(0,T;\widetilde{H}^r(\Omega))$.

To prove \ref{F:Density in dual of X}, note that
\[
C^\infty([0,T];\widetilde{H}^s(\Omega))
\subset
X^1_{s,r}(\Omega_T)
\]
is dense in $L^2(0,T;\widetilde{L}^2(\Omega))$. Hence the first embedding in \eqref{F:eq: dense embedding} is dense. The second embedding is the canonical embedding induced by the $L^2$-duality pairing. Its density follows by duality from the density of the first embedding.

It remains to prove \ref{F:density in X T}. By standard cutoff and mollification arguments in time,
$C_c^\infty([0,T);\widetilde{H}^s(\Omega))$
is dense in $X^1_{s,r,T}(\Omega_T)$; see \cite[Proposition~B.1]{LZ2023unique} or Appendix~\ref{F:appendix: generlaized vector-valued Sobolev}. Since the eigenfunctions form an orthonormal basis of $L^2(\Omega)$
and, after the usual rescaling by the eigenvalues, a Hilbert basis of $\widetilde H^s(\Omega)$, every function in $C_c^\infty([0,T);\widetilde{H}^s(\Omega))$ can be approximated in $H^1(0,T;\widetilde{H}^s(\Omega))$, and hence in $X^1_{s,r,T}(\Omega_T)$, by finite sums of the form \eqref{F:eq: dense set D T}. This proves the claim.
\end{proof}

 \section{Well-posedness of weakly viscous equations}
\label{F:sec: weak finite energy sol}

The main goal of this section is to establish the well-posedness of the weakly viscous nonlocal wave equation
\begin{equation}
\label{F:eq: existence weakly viscous case}
\begin{cases}
( \partial_t^2 +(-\Delta)^r \partial_t +(-\Delta)^s+q) u=f
    &\text{ in }\Omega_T,\\
u=0
    &\text{ in }(\Omega_e)_T,\\
u(0)=u_0,\ \partial_t u(0)=u_1
    &\text{ in }\Omega.
\end{cases}
\end{equation}
This also yields the well-posedness of the corresponding nonhomogeneous exterior value problem
\begin{equation}
\label{F:eq: non-hom weakly viscous case}
\begin{cases}
( \partial_t^2 +(-\Delta)^r \partial_t +(-\Delta)^s+q) u=0
    &\text{ in }\Omega_T,\\
u=\varphi
    &\text{ in }(\Omega_e)_T,\\
u(0)=u_0,\ \partial_t u(0)=u_1
    &\text{ in }\Omega.
\end{cases}
\end{equation}

In contrast to the damped case $r=0$, treated in
\cite{zimmermann2024optimalrungeapproximationdamped}, the term $(-\Delta)^r \partial_t u$ cannot be viewed, a priori, as a lower-order perturbation of the nonlocal wave operator
\[
\partial_t^2+(-\Delta)^s.
\]
A particular feature of the weakly viscous regime $0< r<s$ is the fact that finite energy solutions naturally satisfy
\[
u\in L^2(0,T;\widetilde H^s(\Omega)),
\qquad
\partial_t u\in L^2(0,T;\widetilde H^r(\Omega)).
\]
Thus, the velocity $\partial_t u$ has lower spatial regularity than the displacement $u$, while higher regularity as solutions to nonlocal wave equations without viscous terms. 

To overcome this difficulty, we use a parabolic regularization. More precisely, for $\eps>0$, we consider
\begin{equation}
\label{F:eq: regularized weakly viscous case}
\begin{cases}
( \partial_t^2+\eps(-\Delta)^s\partial_t+(-\Delta)^r \partial_t+(-\Delta)^s+q)u=f
    &\text{ in }\Omega_T,\\
u=0
    &\text{ in }(\Omega_e)_T,\\
u(0)=u_0,\ \partial_t u(0)=u_1
    &\text{ in }\Omega.
\end{cases}
\end{equation}
The additional term $\eps(-\Delta)^s\partial_t u$ gives the improved regularity
\[
\partial_t u_\eps\in L^2(0,T;\widetilde H^s(\Omega)),
\]
which allows us to justify energy identities at the approximate level. We then derive estimates independent of $\eps$ and pass to the limit as $\eps\to0$.

More precisely, in Section~\ref{F:subsec: parabolically regularized}, we recall the abstract viscous evolution framework used to solve \eqref{F:eq: regularized weakly viscous case}. In Section~\ref{F:sec: convergence u eps}, we prove that the solutions $u_\eps$ converge, along a subsequence, to a weak finite energy solution
\begin{equation}
\label{F:eq: reg weak sol}
u\in L^2(0,T;\widetilde H^s(\Omega))
\cap
H^1(0,T;\widetilde H^r(\Omega))
\end{equation}
of \eqref{F:eq: existence weakly viscous case}. Subsequently, in Section~\ref{F:sec: continuity and energy id}, we show that this solution enjoys the stronger regularity
\begin{equation}
\label{F:eq: intro sec existence regularity of sols}
u\in C([0,T];\widetilde H^s(\Omega))
\cap
C^1([0,T];\widetilde L^2(\Omega)).
\end{equation}
We also prove a suitable energy identity; see Proposition~\ref{F:prop: energy identity weakly viscous nonlocal wave}. This identity and the additional time regularity in \eqref{F:eq: intro sec existence regularity of sols} will be important later in the construction of very weak solutions in \cite{zimmermann-paperII}.

Finally, uniqueness of weak finite energy solutions is proved in Section~\ref{F:sec: uniqueness finite energy sol}, and the nonhomogeneous exterior problem \eqref{F:eq: non-hom weakly viscous case} is treated in Section~\ref{F:sec: nonhomgeneous weakly viscous}.

\subsection{Parabolic regularization}
\label{F:subsec: parabolically regularized}

The well-posedness of the regularized problem \eqref{F:eq: regularized weakly viscous case} follows from the abstract theory of viscous evolution equations developed in \cite[Chapter~XVIII]{DautrayLionsVol5}. We recall the relevant framework in a form adapted to the present setting.

Let $V\hookrightarrow H\hookrightarrow V'$ be a Gelfand triple of complex Hilbert spaces, where $V'$ denotes the antidual of $V$. We consider two families of sesquilinear forms
\[
a(t;\cdot,\cdot),\qquad b(t;\cdot,\cdot),
\qquad t\in[0,T],
\]
defined on $V$, and a family of sesquilinear forms
\[
c(t;\cdot,\cdot),\qquad t\in[0,T],
\]
defined on $H$. The corresponding abstract viscous evolution problem is
\begin{equation}
\label{F:eq: abstract evolution problem}
\begin{cases}
\dfrac{d}{dt}c(t;\partial_t u,v)
+b(t;\partial_t u,v)
+a(t;u,v)
=
\langle f(t),v\rangle_{V'\times V}
&\text{in }(0,T),\\
u(0)=u_0,\ \partial_t u(0)=u_1,
\end{cases}
\end{equation}
for all $v\in V$. Here $u_0\in V$, $u_1\in H$, and $f\in L^2(0,T;V')$. We highlight that \eqref{F:eq: abstract evolution problem} has to be understood in the sense of distributions on $(0,T)$.

\begin{definition}
\label{F:def: usual conditions}
We say that a tuple $(V,H,a,b,c)$ satisfies the \emph{usual conditions} if the following assumptions hold.

\begin{enumerate}[{(S}1)]
\item\label{F:S0}
The spaces $V$ and $H$ form a Gelfand triple
\begin{equation}
\label{F:eq: Hilbert triplet condition}
V\hookrightarrow H\hookrightarrow V',
\end{equation}
where $H$ is identified with its antidual via the Riesz representation theorem, and both embeddings are continuous and dense.

\item\label{F:S1}
The family $a(t;\cdot,\cdot)$, $t\in[0,T]$, consists of sesquilinear forms over $V$ and admits a decomposition
\[
a=a_0+a_1.
\]
The principal part $a_0$ satisfies:
\begin{enumerate}[{(A}1)]
\item\label{F:A1}
For all $u,v\in V$, the map $t\mapsto a_0(t;u,v)$ belongs to $C^1([0,T])$, and $a_0(t;\cdot,\cdot)$ is continuous over $V$.
\item\label{F:A2}
The form $a_0$ is Hermitian, that is,
\[
a_0(t;u,v)=\overline{a_0(t;v,u)}
\]
for all $u,v\in V$.
\item\label{F:A3}
The form $a_0$ is coercive over $V$ with respect to $H$, namely,
\[
a_0(t;u,u)\geq \alpha\|u\|_V^2-\lambda\|u\|_H^2
\]
for some $\alpha>0$, $\lambda\in\R$, and all $t\in[0,T]$, $u\in V$.
\end{enumerate}
The lower-order part $a_1$ satisfies:
\begin{enumerate}[(A4)]
\item\label{F:A4}
For all $u,v\in V$, the map $t\mapsto a_1(t;u,v)$ belongs to $C([0,T])$.
\end{enumerate}
\begin{enumerate}[(A5)]
\item\label{F:A5}
There exists $C>0$ such that
\[
|a_1(t;u,v)|
\leq
C\|u\|_V\|v\|_H
\]
for all $u,v\in V$ and $t\in[0,T]$.
\end{enumerate}

\item\label{F:S2}
The family $b(t;\cdot,\cdot)$, $t\in[0,T]$, consists of sesquilinear forms over $V$ and admits a decomposition
\[
b=b_0+b_1.
\]
The forms $b_0,b_1$ satisfy:
\begin{enumerate}[{(B}1)]
\item\label{F:B1}
The family $b_0$ is continuous over $V$, Hermitian, and strongly coercive over $V$, that is,
\[
b_0(t;u,u)\geq \nu\|u\|_V^2
\]
for some $\nu>0$ and all $t\in[0,T]$, $u\in V$.
\item\label{F:B2}
There exists $C>0$ such that
\[
|b_1(t;u,v)|
\leq
C\|u\|_V\|v\|_H
\]
for all $u,v\in V$ and $t\in[0,T]$.
\item\label{F:B3}
For $j=0,1$, the maps $t\mapsto b_j(t;u,v)$ belong to $C([0,T])$ for all $u,v\in V$.
\end{enumerate}

\item\label{F:S3}
The family $c(t;\cdot,\cdot)$, $t\in[0,T]$, consists of sesquilinear forms over $H$ of the form
\[
c(t;u,v)=\langle C(t)u,v\rangle_H,
\]
where $C(t)\in L(H)$. Moreover:
\begin{enumerate}[{(C}1)]
\item\label{F:C1}
The operator $C(t)$ is Hermitian and strongly coercive over $H$, that is,
\[
\langle C(t)u,u\rangle_H\geq \kappa\|u\|_H^2
\]
for some $\kappa>0$ and all $t\in[0,T]$, $u\in H$.
\item\label{F:C2}
For all $u,v\in H$, the map
\[
t\mapsto \langle C(t)u,v\rangle_H
\]
belongs to $C^1([0,T])$.
\end{enumerate}
\end{enumerate}
\end{definition}

\begin{remark}
\label{F:remark: terminology}
We call \eqref{F:eq: abstract evolution problem} an abstract \emph{viscous} evolution equation because the strong coercivity of $b_0$ yields regularization of the velocity field $\partial_t u$, namely
\[
\partial_t u\in L^2(0,T;V).
\]
\end{remark}

\begin{remark}
\label{F:rem: continuity of forms}
By the uniform boundedness principle, condition~\ref{F:A1} implies that $\partial_ta_0(t;\cdot,\cdot)$ is a continuous sesquilinear form over $V$. In particular, there exists $C>0$ such that
\[
|a_0(t;u,v)|+
|\partial_ta_0(t;u,v)|
\leq
C\|u\|_V\|v\|_V
\]
for all $u,v\in V$ and $t\in[0,T]$.
\end{remark}

\begin{remark}
\label{F:remark: inverse of C}
The assumptions on $C(t)$, together with the Lax--Milgram theorem, imply that each $C(t)$ is invertible on $H$. Moreover, $C(t)^{-1}\in L(H)$ is Hermitian and strongly coercive. In addition, the identity
\[
C(t)^{-1}-C(s)^{-1}
=
C(t)^{-1}(C(s)-C(t))C(s)^{-1}
\]
shows that, for all $u,v\in H$,
\[
t\mapsto \langle C(t)^{-1}u,v\rangle_H
\]
belongs to $C^1([0,T])$, and
\[
\frac{d}{dt}C(t)^{-1}
=
-C(t)^{-1}C'(t)C(t)^{-1}.
\]
Consequently,
\[
\begin{split}
|\langle C(t)u,v\rangle_H|
+
|\langle C'(t)u,v\rangle_H|
&\leq C\|u\|_H\|v\|_H,\\
|\langle C(t)^{-1}u,v\rangle_H|
+
|\langle \partial_t C(t)^{-1}u,v\rangle_H|
&\leq C\|u\|_H\|v\|_H .
\end{split}
\]
\end{remark}

We now apply this framework to the regularized weakly viscous equation.

\begin{example}[Gelfand triple]
\label{F:example spaces}
Let $\Omega\subset\R^n$ be a bounded Lipschitz domain and let $s\geq0$. Then
\[
V\vcentcolon=\widetilde H^s(\Omega),
\qquad
H\vcentcolon=\widetilde L^2(\Omega)\vcentcolon =\widetilde H^0(\Omega),\qquad V'=H^{-s}(\Omega)
\]
form a Gelfand triple
\[
V\hookrightarrow H\hookrightarrow V'.
\]
Both embeddings are continuous and dense. Throughout this work, we equip $\widetilde H^s(\Omega)$ with the equivalent norm
\[
    \|u\|_V \vcentcolon = \|u\|_{\widetilde H^s(\Omega)}\vcentcolon=
\|(-\Delta)^{s/2}u\|_{L^2(\R^n)}.
\]
\end{example}

Next, we introduce the class of potentials we are going to work with.

 \begin{definition}
    \label{F:class of potentials}
        For any $s>0$ and open set $\Omega\subset\R^n$, we denote by $\mathcal{Q}_s$ the class of all real-valued functions $q\in L^1_{loc}(\Omega_T)$ satisfying the following conditions:
    \begin{enumerate}[(i)]
        \item $q\in L^{\infty}(0,T;L^p(\Omega))$ for some exponent $2\leq p\leq \infty$ obeying
    \begin{equation}
    \label{F:eq: exponent p}
        \begin{cases}
			n/s\leq p\leq \infty, &\, \text{if }\, 2s< n,\\
			2<p\leq \infty,  &\, \text{if }\, 2s= n,\\
			2\leq p\leq \infty, &\, \text{if }\, 2s > n.
		\end{cases}
    \end{equation}
        \item $q$ is \emph{distributionally continuous in $t$}, meaning that
    \begin{equation}
    \label{F:eq: continuity in t}
        [0,T]\ni t\mapsto \int_{\Omega}q(x,t)\varphi(x)\,dx\in C([0,T])
    \end{equation}
    for all $\varphi\in C_c^{\infty}(\Omega)$.
    \end{enumerate}
    In addition, we write $q\in \mathcal{Q}_s^p$ when $q\in \mathcal{Q}_s$ has integrability exponent $2\leq p\leq \infty$, that is, $q\in L^{\infty}(0,T;L^p(\Omega))$.
    
    Furthermore, we let $\mathcal{Q}_0$ be the collection of all $q\in L^{\infty}(\Omega_T)$, which are distributionally continuous in $t$.
    \end{definition}

    Potentials in the class $\mathcal{Q}_s$ have the following continuity properties.

\begin{lemma}[{\cite[Theorem~3.1]{zimmermann2024-viscous-wave}}]
\label{F:lemma: cond on q}
Let $\Omega\subset\R^n$ be a bounded Lipschitz domain, $T>0$, $s\geq0$, and suppose that $q\in \mathcal{Q}_s^p$. Then
\[
t\mapsto \langle q(t)u,v\rangle_{L^2(\Omega)}
\]
belongs to $C([0,T])$ for all $u,v\in\widetilde{H}^s(\Omega)$. Moreover,
\begin{equation}
\label{F:eq: estimate a1}
\left|
\langle q(t)u,v\rangle_{L^2(\Omega)}
\right|
\lesssim
\|q(t)\|_{L^p(\Omega)}
\|u\|_{\widetilde H^s(\Omega)}
\|v\|_{L^2(\Omega)}
\end{equation}
for all $u,v\in\widetilde H^s(\Omega)$ and a.e. $0\leq t\leq T$. The implicit constant depends only on $n,s,p$, and $\Omega$.
\end{lemma}

\begin{example}[Setup for weakly viscous equations]
\label{F:example sesquilinear forms}
Let $\Omega\subset\R^n$ be a bounded Lipschitz domain, let $0\leq r<s<\infty$ with $s>0$, and let $V,H$ be as in Example~\ref{F:example spaces}. Define
\[
a=a_0+a_1
\]
by
\[
a_0(u,v)
\vcentcolon=
\langle (-\Delta)^{s/2}u,(-\Delta)^{s/2}v\rangle_{L^2(\R^n)}
\]
and
\[
a_1(t;u,v)
\vcentcolon=
\langle q(t)u,v\rangle_{L^2(\Omega)}.
\]
Then $a_0$ satisfies \ref{F:A1}--\ref{F:A3}, and
\[
a_0(u,u)=\|u\|_V^2.
\]
Moreover, by Lemma~\ref{F:lemma: cond on q}, the family $a_1$ satisfies \ref{F:A4}--\ref{F:A5}.

For $\eps>0$, define
\[
b^{r,\eps}=b_0^{r,\eps}+b_1^{r,\eps}
\]
as follows. If $0\leq r\leq s/2$, set
\[
b_0^{r,\eps}(u,v)
\vcentcolon=
\eps a_0(u,v),
\qquad
b_1^{r,\eps}(u,v)
\vcentcolon=
\langle (-\Delta)^r u,v\rangle_{L^2(\Omega)}.
\]
If $s/2<r<s$, set
\[
b_0^{r,\eps}(u,v)
\vcentcolon=
\eps a_0(u,v)
+
\langle (-\Delta)^{r/2}u,(-\Delta)^{r/2}v\rangle_{L^2(\R^n)},
\qquad
b_1^{r,\eps}\equiv0.
\]
In both cases the forms $b_0^{r,\eps}$ and $b_1^{r,\eps}$ satisfy \ref{F:B1}--\ref{F:B3}. Indeed, in the case $0\leq r\leq s/2$, the embedding $H^s(\R^n)\hookrightarrow H^{2r}(\R^n)$ gives
\[
|b_1^{r,\eps}(u,v)|
\leq
\|(-\Delta)^r u\|_{L^2(\R^n)}
\|v\|_{L^2(\Omega)}
\leq
C\|u\|_V\|v\|_H.
\]
In the case $s/2<r<s$, the assertion follows from the embedding $H^s(\R^n)\hookrightarrow H^r(\R^n)$.

Finally, let
\[
c(u,v)\vcentcolon=\langle u,v\rangle_H.
\]
Then the tuple $(V,H,a,b^{r,\eps},c)$ satisfies the usual conditions.
\end{example}

\begin{remark}[Damping term]
\label{F:rem: damping term}
One may include an additional damping term $\gamma\partial_t u$ in \eqref{F:eq: regularized weakly viscous case} by replacing $b_1^{r,\eps}$ by
\[
b_1^{r,\eps}(t;u,v)
=
\langle (-\Delta)^r u,v\rangle_{L^2(\Omega)}
+
\langle \gamma(t)u,v\rangle_{L^2(\Omega)}
\]
when $0\leq r\leq s/2$, and by
\[
b_1^{r,\eps}(t;u,v)
=
\langle \gamma(t)u,v\rangle_{L^2(\Omega)}
\]
when $s/2<r<s$. The conditions \ref{F:B2}--\ref{F:B3} follow, for instance, if $\gamma\in\mathcal Q_s$. Additional positivity or regularity assumptions on $\gamma$ may be needed to obtain uniform energy estimates as $\eps\to0$.
\end{remark}

We shall also use the following solution spaces associated with the abstract viscous and non-viscous evolution problems.

\begin{definition}
\label{F:def: special function spaces}
Let $V\hookrightarrow H\hookrightarrow V'$ be a Gelfand triple and let $c(t;\cdot,\cdot)$ be the sesquilinear forms on $H$ induced by the operators $C(t)\in L(H)$, $t\in[0,T]$, satisfying condition~\ref{F:S3}. Then we define
\[
W_c(0,T;V)
\vcentcolon=
\{u\in L^2(0,T;V)\,;\,\partial_t(Cu)\in L^2(0,T;V')\},
\]
\[
\widetilde W_c(0,T;V)
\vcentcolon=
\{u\in L^2(0,T;V)\,;\,\partial_t u\in W_c(0,T;V)\},
\]
and
\[
\widetilde{\widetilde W}_c(0,T;V)
\vcentcolon=
\{u\in L^2(0,T;V)\,;\,
\partial_t u\in L^2(0,T;H),\
\partial_t(C\partial_t u)\in L^2(0,T;V')\}.
\]
These spaces are equipped with their natural graph norms. If $C(t)=\mathrm{id}_H$, we omit the subscript $c$.
\end{definition}

\begin{lemma}[{\cite[Chapter~XVIII, \S~5]{DautrayLionsVol5}}]
\label{F:lemma: embeddings}
Assume that the hypotheses of Definition~\ref{F:def: special function spaces} are satisfied. Then:
\begin{enumerate}[(i)]
\item
The embeddings
\[
W_c(0,T;V)\hookrightarrow C([0,T];H),
\qquad
\widetilde W_c(0,T;V)\hookrightarrow C([0,T];V)
\]
are continuous.
\item
The space $C_c^\infty([0,T];V)$ is dense in both $W_c(0,T;V)$ and $\widetilde W_c(0,T;V)$.
\item
For all $u\in\widetilde W_c(0,T;V)$,
\begin{equation}
\label{F:eq: integration by parts with c 1}
\begin{split}
&2\re\int_0^T
\langle \partial_t(C\partial_t u),\partial_t u\rangle_{V'\times V}\,dt
-
\int_0^T
\langle C'\partial_t u,\partial_t u\rangle_H\,dt
\\
&=
c(T;\partial_t u(T),\partial_t u(T))
-
c(0;\partial_t u(0),\partial_t u(0)).
\end{split}
\end{equation}
\end{enumerate}
\end{lemma}

We can now state the abstract well-posedness result.

\begin{theorem}[Well-posedness of abstract viscous evolution equations]
\label{F:thm: abstract well-posedness}
Assume that $(V,H,a,b,c)$ satisfies the usual conditions. Then, for every
\[
 (u_0,u_1,f)\in V\times H\times L^2(0,T;V'),
\]
there exists a unique solution $u\in \widetilde W_c(0,T;V)$ of \eqref{F:eq: abstract evolution problem}. Moreover, $u$ satisfies the energy identity
\begin{equation}
\label{F:eq: abstract energy identity}
\begin{split}
\Psi(t)
&=
\Psi(0)
+
\int_0^t a_0'(\tau;u(\tau),u(\tau))\,d\tau
\\
&\quad
-2\re\int_0^t a_1(\tau;u(\tau),\partial_t u(\tau))\,d\tau
-2\re\int_0^t b_1(\tau;\partial_t u(\tau),\partial_t u(\tau))\,d\tau
\\
&\quad
-\int_0^t
\langle C'(\tau)\partial_t u(\tau),\partial_t u(\tau)\rangle_H\,d\tau
+
2\re\int_0^t
\langle f(\tau),\partial_t u(\tau)\rangle_{V'\times V}\,d\tau,
\end{split}
\end{equation}
where
\[
\Psi(t)
\vcentcolon=
c(t;\partial_t u(t),\partial_t u(t))
+
a_0(t;u(t),u(t))
+
2\int_0^t b_0(\tau;\partial_t u(\tau),\partial_t u(\tau))\,d\tau .
\]
\end{theorem}

\begin{proof}
This follows from \cite[Chapter~XVIII, \S~5, Theorem~1 and Remark~4]{DautrayLionsVol5}. The identity \eqref{F:eq: abstract energy identity} follows from \cite[Chapter~XVIII, \S~5, eq.~(5.83)]{DautrayLionsVol5}.
\end{proof}

As an immediate consequence, we obtain the well-posedness of the parabolically regularized equation.

\begin{corollary}[Well-posedness of the parabolically regularized equation]
\label{F:cor: regularized weakly viscous nonlocal wave}
Let $\Omega\subset\R^n$ be a bounded Lipschitz domain, $T>0$, $\eps>0$, and $0\leq r<s<\infty$, $s>0$. Suppose that $q\in\mathcal Q_s$. Then, for every
\[
(u_0,u_1,f)
\in
\widetilde H^s(\Omega)
\times
\widetilde L^2(\Omega)
\times
L^2(0,T;H^{-s}(\Omega)),
\]
there exists a unique solution $u_\eps\in \widetilde W(0,T;\widetilde H^s(\Omega))$ of \eqref{F:eq: regularized weakly viscous case}. Moreover, $u_\eps$ satisfies
\begin{equation}
\label{F:eq: energy identity regularized viscous}
\begin{split}
&\|\partial_t u_\eps(t)\|_{L^2(\Omega)}^2
+
\|(-\Delta)^{s/2}u_\eps(t)\|_{L^2(\R^n)}^2+2\int_0^t
\|(-\Delta)^{r/2} \partial_t u_\eps(\tau)\|_{L^2(\R^n)}^2\,d\tau
\\
&\quad
+
2\eps\int_0^t
\|(-\Delta)^{s/2}\partial_t u_\eps(\tau)\|_{L^2(\R^n)}^2\,d\tau
\\
&=
\|u_1\|_{L^2(\Omega)}^2
+
\|(-\Delta)^{s/2}u_0\|_{L^2(\R^n)}^2-2\re\int_0^t
\langle q(\tau)u_\eps(\tau),\partial_t u_\eps(\tau)\rangle_{L^2(\Omega)}\,d\tau
\\
&\quad
+
2\re\int_0^t
\langle f(\tau),\partial_t u_\eps(\tau)\rangle_{H^{-s}(\Omega)\times\widetilde H^s(\Omega)}\,d\tau .
\end{split}
\end{equation}
\end{corollary}

\begin{proof}
This follows from Theorem~\ref{F:thm: abstract well-posedness} applied to the setting of Example~\ref{F:example sesquilinear forms}. Since
\[
    u_\eps\in \widetilde W(0,T;\widetilde H^s(\Omega))= H^1(0,T;\widetilde H^s(\Omega))\cap H^2(0,T;H^{-s}(\Omega)),
\]
by classical embedding theorems, all terms in \eqref{F:eq: energy identity regularized viscous} are well-defined.
\end{proof}

\subsection{Convergence of the solutions \texorpdfstring{$u_\eps$}{u eps}}
\label{F:sec: convergence u eps}

We now pass to the limit $\eps\to 0$. The key point is to derive estimates that are uniform in $\eps$. We will use repeatedly Gronwall's inequality in the following form (see Section~\ref{F:appendix: Proof of Gronwall}).

\begin{lemma}[Gronwall's inequality]
\label{F:lemma: Gronwall}
Let $T>0$. Assume that $\varphi,\alpha\in L^{\infty}(0,T)$ and $\mu\in L^1(0,T)$ are nonnegative functions such that $\alpha$ is nondecreasing and
\begin{equation}
\label{F:eq: gronwall assumpt}
\varphi(t)
\leq
\alpha(t)
+
\int_0^t \mu(\tau)\varphi(\tau)\,d\tau
\end{equation}
for a.e. $0<t<T$. Then
\[
\varphi(t)
\leq
\alpha(t)
\exp\left(\int_0^t\mu(\tau)\,d\tau\right)
\]
for a.e. $0<t<T$.
\end{lemma}

\begin{lemma}[Uniform estimate]
\label{F:lemma: uniform bound regularized sols}
Let $\Omega\subset\R^n$ be a bounded Lipschitz domain, $T>0$, and $0\leq r<s<\infty$. Suppose that $q\in \mathcal{Q}_s^p$,
\[
(u_0,u_1,f)\in\widetilde H^s(\Omega)\times\widetilde L^2(\Omega)\times L^2(0,T;H^{-r}(\Omega)).
\]
Then the solutions $u_\eps$, $\eps>0$, of \eqref{F:eq: regularized weakly viscous case} satisfy
\begin{equation}
\label{F:eq: uniform bound on regularized sols}
\begin{split}
&\|\partial_t u_\eps(t)\|_{L^2(\Omega)}^2
+
\|u_\eps(t)\|_{\widetilde H^s(\Omega)}^2
+
\int_0^t
\|\partial_t u_\eps(\tau)\|_{\widetilde H^r(\Omega)}^2\,d\tau
+
\eps\int_0^t
\|\partial_t u_\eps(\tau)\|_{\widetilde H^s(\Omega)}^2\,d\tau
\\
&\leq
C
\exp\left(
Ct(1+\|q\|_{L^\infty(0,t;L^p(\Omega))}^2)
\right)
\left(
\|u_1\|_{L^2(\Omega)}^2
+
\|u_0\|_{\widetilde H^s(\Omega)}^2
+
\|f\|_{L^2(0,t;H^{-r}(\Omega))}^2
\right)
\end{split}
\end{equation}
for all $0\leq t\leq T$, where $C>0$ is independent of $\eps$.
\end{lemma}

\begin{proof}
Using \eqref{F:eq: energy identity regularized viscous}, Lemma~\ref{F:lemma: cond on q}, and Young's inequality, we obtain
\[
\begin{split}
&\|\partial_t u_\eps(t)\|_{L^2(\Omega)}^2
+
\|u_\eps(t)\|_{\widetilde H^s(\Omega)}^2
+
2\int_0^t
\|\partial_t u_\eps(\tau)\|_{\widetilde H^r(\Omega)}^2\,d\tau
\\
&\quad
+
2\eps\int_0^t
\|\partial_t u_\eps(\tau)\|_{\widetilde H^s(\Omega)}^2\,d\tau
\\
&\leq
\|u_1\|_{L^2(\Omega)}^2
+
\|u_0\|_{\widetilde H^s(\Omega)}^2
\\
&\quad
+
C\int_0^t
\bigl(1+\|q(\tau)\|_{L^p(\Omega)}^2\bigr)
\left(
\|\partial_t u_\eps(\tau)\|_{L^2(\Omega)}^2
+
\|u_\eps(\tau)\|_{\widetilde H^s(\Omega)}^2
\right)\,d\tau
\\
&\quad
+
C\int_0^t
\|f(\tau)\|_{H^{-r}(\Omega)}^2\,d\tau
+
\int_0^t
\|\partial_t u_\eps(\tau)\|_{\widetilde H^r(\Omega)}^2\,d\tau .
\end{split}
\]
We absorb the last term into the left-hand side.

Set
\[
\begin{split}
\Phi_\eps(t)
\vcentcolon={}&
\|\partial_t u_\eps(t)\|_{L^2(\Omega)}^2
+
\|u_\eps(t)\|_{\widetilde H^s(\Omega)}^2
\\
&+
\int_0^t
\|\partial_t u_\eps(\tau)\|_{\widetilde H^r(\Omega)}^2\,d\tau
+
\eps\int_0^t
\|\partial_t u_\eps(\tau)\|_{\widetilde H^s(\Omega)}^2\,d\tau .
\end{split}
\]
Then we deduce that
\[
\Phi_\eps(t)
\leq
C\left(
\Phi_\eps(0)
+
\|f\|_{L^2(0,t;H^{-r}(\Omega))}^2
\right)
+
C\int_0^t
\bigl(1+\|q(\tau)\|_{L^p(\Omega)}^2\bigr)
\Phi_\eps(\tau)\,d\tau 
\]
and Gronwall's inequality yields \eqref{F:eq: uniform bound on regularized sols}.
\end{proof}

\begin{proposition}[Limit as $\eps\to0$]
\label{F:prop: limit function}
Assume that the hypotheses of Lemma~\ref{F:lemma: uniform bound regularized sols} hold. Then there exists
\[
u\in L^\infty(0,T;\widetilde H^s(\Omega))
\cap
H^1(0,T;\widetilde H^r(\Omega))
\]
with
\[
\partial_t u\in L^\infty(0,T;\widetilde L^2(\Omega))
\]
such that, along a subsequence,
\[
\begin{split}
    u_\eps\weakstar u
\quad&\text{in }L^\infty(0,T;\widetilde H^s(\Omega)),\\
\partial_t u_\eps\weakstar \partial_t u
\quad&\text{in }L^\infty(0,T;\widetilde L^2(\Omega)),\\
\partial_t u_\eps\rightharpoonup \partial_t u
\quad&\text{in }L^2(0,T;\widetilde H^r(\Omega)).
\end{split}
\]
Moreover, $u$ solves \eqref{F:eq: existence weakly viscous case} in the weak sense, namely,
\begin{equation}
\label{F:eq: calc for distr sol}
\begin{split}
&\int_0^T
\left(
\langle u,\partial_t^2 \varphi\rangle_{L^2(\Omega)}
+
\langle (-\Delta)^{r/2} \partial_t u,(-\Delta)^{r/2}\varphi\rangle_{L^2(\R^n)}
\right)\,dt
\\
&\quad
+
\int_0^T
\left(
\langle (-\Delta)^{s/2}u,(-\Delta)^{s/2}\varphi\rangle_{L^2(\R^n)}
+
\langle q u,\varphi\rangle_{L^2(\Omega)}
\right)\,dt
\\
&=
\langle u_1,\varphi(0)\rangle_{L^2(\Omega)}
-
\langle u_0,\partial_t\varphi(0)\rangle_{L^2(\Omega)}
+
\int_0^T
\langle f,\varphi\rangle_{H^{-r}(\Omega)\times\widetilde H^r (\Omega)}\,dt
\end{split}
\end{equation}
for any $\varphi\in C_c^{\infty}(\Omega\times [0,T))$. Furthermore, $u$ satisfies
\[
u(0)=u_0
\text{ in }\widetilde H^r(\Omega),
\]
and obeys the estimate
\begin{equation}
\label{F:eq: energy estimate weakly viscous}
\begin{split}
&\|\partial_t u(t)\|_{L^2(\Omega)}^2
+
\|u(t)\|_{\widetilde H^s(\Omega)}^2
+
\int_0^t
\|\partial_t u(\tau)\|_{\widetilde H^r(\Omega)}^2\,d\tau
\\
&\leq
C
\exp\left(
C t(1+\|q\|_{L^\infty(0,t;L^p(\Omega))}^2)
\right)
\left(
\|u_1\|_{L^2(\Omega)}^2
+
\|u_0\|_{\widetilde H^s(\Omega)}^2
+
\|f\|_{L^2(0,t;H^{-r}(\Omega))}^2
\right)
\end{split}
\end{equation}
for a.e. $0\leq t\leq T$.
\end{proposition}

\begin{proof}
The uniform estimate \eqref{F:eq: uniform bound on regularized sols} implies that
\[
(u_\eps)_{\eps>0}
\text{ is bounded in }L^\infty(0,T;\widetilde H^s(\Omega)),
\]
\[
(\partial_t u_\eps)_{\eps>0}
\text{ is bounded in }
L^\infty(0,T;\widetilde L^2(\Omega))
\cap
L^2(0,T;\widetilde H^r(\Omega)),
\]
and
\[
(\eps^{1/2}\partial_t u_\eps)_\eps
\text{ is bounded in }
L^2(0,T;\widetilde H^s(\Omega)).
\]
Hence, after passing to a subsequence, we obtain the stated weak and weak-$\ast$ convergences. Since
\[
u_\eps\rightharpoonup u
\quad\text{in }H^1(0,T;\widetilde H^r(\Omega)),
\]
the trace operator
\[
H^1(0,T;\widetilde H^r(\Omega))
\to
\widetilde H^r(\Omega),\,u\mapsto u(0),
\]
is continuous, and therefore
\[
u(0)=u_0
\text{ in }\widetilde H^r(\Omega).
\]

Let $\varphi\in C_c^\infty(\Omega\times[0,T))$. Since $u_\eps$, $\eps>0$, solves \eqref{F:eq: regularized weakly viscous case}, integration by parts in time gives
\[
\begin{split}
&\int_0^T
\left(
\langle u_\eps,\partial_t^2\varphi\rangle_{L^2(\Omega)}
+
\langle (-\Delta)^{r/2} \partial_t u_\eps,(-\Delta)^{r/2}\varphi\rangle_{L^2(\R^n)}
\right)\,dt
\\
&\quad
+
\int_0^T
\left(
\langle (-\Delta)^{s/2}u_\eps,(-\Delta)^{s/2}\varphi\rangle_{L^2(\R^n)}
+
\langle q u_\eps,\varphi\rangle_{L^2(\Omega)}
\right)\,dt
\\
&\quad
+
\eps\int_0^T
\langle (-\Delta)^{s/2}\partial_t u_\eps,
(-\Delta)^{s/2}\varphi\rangle_{L^2(\R^n)}\,dt
\\
&=
\langle u_1,\varphi(0)\rangle_{L^2(\Omega)}
-
\langle u_0,\partial_t\varphi(0)\rangle_{L^2(\Omega)}
+
\int_0^T
\langle f,\varphi\rangle_{H^{-r}(\Omega)\times\widetilde H^r(\Omega)}\,dt .
\end{split}
\]
The last term on the left-hand side tends to zero because
\[
\eps^{1/2}\partial_t u_\eps
\]
is bounded in $L^2(0,T;\widetilde H^s(\Omega))$. Passing to the limit in the remaining terms yields \eqref{F:eq: calc for distr sol}. Note that the convergence of the potential term follows from Lemma~\ref{F:lemma: cond on q} and the weak-$\ast$ convergence of $u_\eps$ in $L^\infty(0,T;\widetilde H^s(\Omega))$.

Thus $u$ is a weak solution of \eqref{F:eq: existence weakly viscous case}. Finally, the estimate \eqref{F:eq: energy estimate weakly viscous} follows from \eqref{F:eq: uniform bound on regularized sols} by weak lower semicontinuity.
\end{proof}

\begin{definition}
\label{F:def: def of solutions}
Assume that the hypotheses of Lemma~\ref{F:lemma: uniform bound regularized sols} hold. We call $u$ a \emph{weak finite energy solution} of \eqref{F:eq: existence weakly viscous case} if
\[
u\in L^2(0,T;\widetilde H^s(\Omega))
\cap
H^1(0,T;\widetilde H^r(\Omega)),
\]
$u$ weakly solves \eqref{F:eq: existence weakly viscous case}, meaning that it satisfies \eqref{F:eq: calc for distr sol} for any $\varphi\in C_c^{\infty}(\Omega\times [0,T))$, and
\[
u(0)=u_0
\text{ in }\widetilde H^r(\Omega).
\]
We call $u$ a \emph{finite energy solution} if, in addition,
\[
u\in C([0,T];\widetilde H^s(\Omega))
\cap
C^1([0,T];\widetilde L^2(\Omega))
\]
and
\[
(u(0),\partial_t u(0))
=
(u_0,u_1)
\text{ in }
\widetilde H^s(\Omega)\times\widetilde L^2(\Omega).
\]
\end{definition}

For later use, we record the equivalent weak formulation. If $u$ is a weak finite energy solution, then
\begin{equation}
\label{F:eq: weak formualation}
\begin{split}
&\int_0^T
\left(
-\langle \partial_t u,\partial_t\varphi\rangle_{L^2(\Omega)}
+
\langle (-\Delta)^{r/2} \partial_t u,(-\Delta)^{r/2}\varphi\rangle_{L^2(\R^n)}
\right)\,dt
\\
&\quad
+
\int_0^T
\left(
\langle (-\Delta)^{s/2}u,(-\Delta)^{s/2}\varphi\rangle_{L^2(\R^n)}
+
\langle qu,\varphi\rangle_{L^2(\Omega)}
\right)\,dt
\\
&=
\langle u_1,\varphi(0)\rangle_{L^2(\Omega)}
+
\int_0^T
\langle f,\varphi\rangle_{H^{-r} (\Omega)\times\widetilde H^r(\Omega)}\,dt
\end{split}
\end{equation}
for all $\varphi\in C_c^\infty(\Omega\times[0,T))$.

    \subsection{Continuity and energy identity}
    \label{F:sec: continuity and energy id}

    In this section, we aim to establish the following two properties:
    \begin{enumerate}[(P1)]
        \item\label{F:P1} The energy estimate \eqref{F:eq: energy estimate weakly viscous} can be sharpened into an energy identity.
    \end{enumerate}
    \begin{enumerate}[(P2)]
        \item\label{F:P2} The constructed solution of \eqref{F:eq: existence weakly viscous case}  satisfies
        \begin{equation}
        \label{F:eq: continuity property}
            u\in C([0,T];\widetilde{H}^s(\Omega))\cap C^1([0,T];\widetilde{L}^2(\Omega))
        \end{equation}
        and the initial conditions are attained pointwise in $\widetilde{H}^s(\Omega)\times\widetilde{L}^2(\Omega)$.
    \end{enumerate}
    Let us note that the energy identity does not directly follow from the energy identity for $u_\eps$ as at this point we only know that $\eps\int_{\R^n}|(-\Delta)^{s/2} \partial_t u_\eps|^2\,dx$ is bounded. But, in principle, the desired energy identity can be derived, as usual, by testing the PDE for $u$ with $\partial_t u$ and performing an integration by parts. Unfortunately, $\partial_t u$ is not regular enough to justify the integration by parts formula (see Lemma~\ref{F:lemma: embeddings}). Therefore, we are forced to follow \cite{lions2012nonVol1} in which this regularity issue is circumvented by performing a double regularization.

    The energy identity turns out to be crucial to establish \eqref{F:eq: continuity property} and to develop an existence theory for nonlinear weakly viscous nonlocal wave equations.

    To carry out the above plan, we need to introduce the space of scalar continuous functions, which are defined similarly to distributionally continuous functions. If $X$ is any normed space with antidual space $X'$ and $T>0$, then we define the space of \emph{scalar continuous ($X$-valued) functions on $[0,T]$} by
    \begin{equation}
    \label{F:eq: scalar continuous functions}
        C_w([0,T];X)=\{u\in L^{\infty}(0,T;X)\,;\,t\mapsto \langle\ell,u(t)\rangle\in C([0,T]),\,\forall\ell\in X'\}.
    \end{equation}
    In contrast to the reference \cite{lions2012nonVol1}, we denote by $C_w([0,T];X)$ the space of scalar continuous functions, instead of using the notation $C_s([0,T];X)$, in order to avoid confusion with the fractional index $s$. We may notice that the scalar continuous functions are defined in complete analogy to the notion of weak convergence of sequences, which explains our choice of notation.
    
    Using classical embedding theorems and Proposition \ref{F:prop: limit function}, we see that the constructed solution $u$ of \eqref{F:eq: existence weakly viscous case} satisfies
    \begin{equation}
    \label{F:eq: continuity of sols}
            u\in C([0,T];\widetilde{L}^2(\Omega))\text{ and }\partial_t u\in C([0,T];H^{-s}(\Omega)).
    \end{equation}
    Next, let us recall the following lemma, whose proof follows by mollification and weak compactness.
    \begin{lemma}[{\cite[Chapter 3, Lemma 8.1]{lions2012nonVol1}}]
    \label{F:lemma: improved continuity}
        Let $T>0$ and suppose that $X,Y$ are two Banach spaces such that $X$ is reflexive and $X\hookrightarrow Y$. Then there holds
        \[
            C_w([0,T];X)=L^{\infty}(0,T;X)\cap C_w([0,T];Y).
        \]
    \end{lemma}
    Thus, combined with \eqref{F:eq: continuity of sols}, we infer that the solution $u$ satisfies
    \begin{equation}
    \label{F:eq: scalar continuity 2}
            u\in C_w([0,T];\widetilde{H}^s(\Omega))\text{ and }\partial_t u\in C_w([0,T];\widetilde{L}^2(\Omega)).
    \end{equation}
    So, to prove property \ref{F:P2} it remains to show that $u$ and $\partial_t u$ are in fact continuous. 

    Before proceeding, we show that our solution $u$ attains the initial data $(u_0,u_1)$ in $\widetilde{H}^s(\Omega)\times\widetilde{L}^2(\Omega)$.

    \begin{lemma}
    \label{F:lemma: u1 is attained strongly}
        Assume the conditions of Lemma \ref{F:lemma: uniform bound regularized sols} hold. If $u$ is a weak finite energy solution of \eqref{F:eq: existence weakly viscous case} with
        \[
            u\in C_w([0,T];\widetilde{H}^s(\Omega))\text{ and }\partial_t u\in C_w([0,T];\widetilde{L}^2(\Omega))\cap L^2(0,T;\widetilde{H}^r(\Omega)),
        \]
        then 
        \[
            (u(0),\partial_t u(0))=(u_0,u_1)\text{ in }\widetilde{H}^s(\Omega)\times \widetilde{L}^2(\Omega).
        \]
    \end{lemma}

    \begin{proof}
        For any $0<h\ll T$, we introduce the functions
        \begin{equation}
        \label{F:eq: eta h}
            \eta_h\colon \R\to \R\text{ with }\eta_h(t)\vcentcolon =\begin{cases}
                1-t/h,&\text{ for }0\leq t\leq h\\
                0,&\text{ otherwise}.
            \end{cases}
        \end{equation}
        By approximation, we can use $\varphi_h(x,t)=\eta_h(t)\psi(x)$, with $\psi\in \widetilde{H}^s(\Omega)$, as test functions in \eqref{F:eq: weak formualation}, which ensures that\footnote{For a measure $\mu$ on the set $X$, a $\mu-$measurable subset $E\subset X$, and $f\in L^1(E,\mu)$, we write 
       \[
            \fint_{E} f\,d\mu\vcentcolon = \frac{1}{\mu(E)}\int_E f \,d\mu.
       \]}
        \[
        \begin{split}
            &\fint_0^h\langle \partial_t u,\psi\rangle_{L^2(\Omega)}\,dt +\int_0^T\langle (-\Delta)^{r/2} \partial_t u,(-\Delta)^{r/2}\psi\rangle_{L^2(\R^n)}\eta_h\,dt\\
            &+\int_0^T\left(\langle (-\Delta)^{s/2}u,(-\Delta)^{s/2}\psi\rangle_{L^2(\R^n)}+\langle q u,\psi\rangle_{L^2(\Omega)}\right)\eta_h\,dt\\
            &= \langle u_1,\psi\rangle_{L^2(\Omega)}+\int_0^T\langle f,\psi\rangle_{L^2(\Omega)}\eta_h\,dt.
        \end{split}
        \]
        By dominated convergence, $\partial_t u\in C_w([0,T];\widetilde{L}^2(\Omega))$, and Lebesgue's differentiation theorem, we get in the limit $h\to 0$ the identity
        \[
            \langle \partial_t u(0),\psi\rangle_{L^2(\Omega)}=\langle u_1,\psi\rangle_{L^2(\Omega)}.
        \]
        Hence, we have shown that $\partial_t u(0)=u_1$ holds in $\widetilde{L}^2(\Omega)$. Next, we prove that $u(0)=u_0$ in $\widetilde{H}^s(\Omega)$. On the one hand, $u\in C_w([0,T];\widetilde{H}^s(\Omega))$ implies that $u(0)\in \widetilde{H}^s(\Omega)$. On the other hand, Proposition~\ref{F:prop: limit function} yields $u(0)=u_0$ in $\widetilde{H}^r(\Omega)$. Thus, $\widetilde{H}^s(\Omega)\hookrightarrow \widetilde{H}^r(\Omega)$ ensures that the equality holds in $\widetilde{H}^s(\Omega)$ as well. Therefore, we can conclude the proof.
    \end{proof}

    Our next result, Proposition~\ref{F:prop: energy identity weakly viscous nonlocal wave}, establishes two different energy identities. Since the proof is rather technical and lengthy, we present it in Appendix~\ref{F:sec: proof of energy identity}. The first formula holds for potentials satisfying the regularity conditions of Lemma~\ref{F:lemma: cond on q}. The second energy identity assumes that the potential $q$ is of class $C^1([0,T];L^p(\Omega))$, but gives further insights into the energy production processes of the weakly viscous wave. More concretely, if we denote the (modified) energy of the wave at time $t$ by
    \[
        \mathcal{E}(t)\vcentcolon =\|\partial_t u(t)\|^2_{L^2(\Omega)}+\|(-\Delta)^{s/2}u(t)\|^2_{L^2(\R^n)}+\langle q(t)u(t),u(t)\rangle_{L^2(\Omega)},
    \]
    then the second energy identity below shows that 
    \[
        \Delta\mathcal{E}(t)\vcentcolon = \mathcal{E}(t)-\mathcal{E}(0)=\int_{\Omega_t}e_q(x,t)\,dxdt-2\int_{\R^n_t}e_r(x,t)\,dxdt,
    \]
    where
    \[
        e_q(x,t)= \partial_t q(x,t)|u(x,t)|^2\quad\text{and}\quad e_r(x,t)=|(-\Delta)^{r/2}\partial_t u(x,t)|^2.
    \]
    Hence, the time-dependence of the potential term may lead to a local energy production whereas the viscosity term induces a (non-) local energy loss.

    \begin{proposition}[Energy identity]
    \label{F:prop: energy identity weakly viscous nonlocal wave}
        Assume the conditions of Lemma \ref{F:lemma: uniform bound regularized sols} hold. Let $u$ be any weak finite energy solution of \eqref{F:eq: existence weakly viscous case} with the continuity properties
        \[
            u\in C_w([0,T];\widetilde{H}^s(\Omega))\text{ and }\partial_t u\in C_w([0,T];\widetilde{L}^2(\Omega)).
        \]
        Then, $u$ satisfies the \emph{energy identity}\footnote{Here and throughout this work, we occasionally drop the subscripts in the antiduality pairings for the sake of readability when they are clear from the context.}
         \begin{equation}
        \label{F:eq: energy identity weakly viscous nonlocal wave eq continuous potential}
               \begin{split}
             &\|\partial_t u(t)\|_{L^2(\Omega)}^2+\|(-\Delta)^{s/2}u(t)\|_{L^2(\R^n)}^2+2 \int_0^t \|
            (-\Delta)^{r/2} \partial_t u\|_{L^2(\R^n)}^2\,d\tau\\
             &=\|u_1\|_{L^2(\Omega)}^2+\|(-\Delta)^{s/2}u_0\|_{L^2(\R^n)}^2 -2 \re \int_0^t \langle qu,\partial_t u\rangle_{L^2(\Omega)}\,d\tau\\
            &\quad +2\re \int_0^t\langle f,\partial_t u\rangle\,d\tau
        \end{split}
        \end{equation}
        for all $0\leq t\leq T$. Furthermore, if the potential $q$ belongs to $C^1([0,T];L^p(\Omega))$, then we have
        \begin{equation}
        \label{F:eq: energy identity weakly viscous nonlocal wave eq}
             \begin{split}
             &\|\partial_t u(t)\|^2_{L^2(\Omega)}+\|(-\Delta)^{s/2}u(t)\|^2_{L^2(\R^n)}+2\int_{0}^{t}\| (-\Delta)^{r/2} \partial_t u\|_{L^2(\R^n)}^2\,d\tau\\
             &+\langle q(t)u(t),u(t)\rangle_{L^2(\Omega)}-\int_{0}^{t}\langle (\partial_tq)u, u\rangle_{L^2(\Omega)}\,d\tau\\
            &=\|u_1\|^2_{L^2(\Omega)}+\|(-\Delta)^{s/2}u_0\|^2_{L^2(\R^n)}+\langle q(0)u_0,u_0\rangle_{L^2(\Omega)}+ 2\re \int_{0}^{t}\langle   f ,\partial_t u\rangle\,d\tau
        \end{split}
        \end{equation}
        for all $0\leq t\leq T$.
    \end{proposition}

 As mentioned above, the main obstruction in proving Proposition~\ref{F:prop: energy identity weakly viscous nonlocal wave} is the limited regularity of the constructed solutions, which prevents a direct justification of the required integrations by parts. To overcome this difficulty, the proof relies on a double approximation argument. More precisely, we first extend the solution $u$ to $\R$, then localize it in time by multiplying it with a cutoff function $\eta_\lambda$, and subsequently regularize the resulting function by convolution with a standard one-dimensional mollifier $\rho_\delta$. The energy identity is then established for the smooth approximation
\[
    \rho_\delta \ast (\eta_\lambda u),
\]
which is compactly supported in time. Finally, by first passing to the limit $\lambda\to0$ and afterwards $\delta\to0$, we recover the energy identities stated in Proposition~\ref{F:prop: energy identity weakly viscous nonlocal wave}.

    As an application of the energy identity, we can now establish the existence of finite energy solutions.
    \begin{theorem}
    \label{F:prop: existence of weak solutions}
        Let $\Omega\subset\R^n$ be a bounded Lipschitz domain, $T>0$, and $0\leq r<s<\infty$. Furthermore, assume that $q\in\mathcal{Q}^p_s$, $(u_0,u_1)\in \widetilde{H}^s(\Omega)\times \widetilde{L}^2(\Omega)$, and $f\in L^2(0,T;H^{-r}(\Omega))$. Then problem \eqref{F:eq: existence weakly viscous case} has a finite energy solution satisfying the energy identity \eqref{F:eq: energy identity weakly viscous nonlocal wave eq continuous potential}.
    \end{theorem}

    \begin{proof}
        Note that the existence of a weak finite energy solution $u$ follows from Proposition \ref{F:prop: limit function}. Furthermore, by the arguments at the beginning of Section \ref{F:sec: continuity and energy id} we know that $u$ has the following continuity properties
        \[
            u\in C_w([0,T];\widetilde{H}^s(\Omega))\text{ and }\partial_t u\in C_w([0,T];\widetilde{L}^2(\Omega))\cap L^2(0,T;\widetilde{H}^r(\Omega))
        \]
        and it attains the initial data $(u_0,u_1)$ in $\widetilde H^s(\Omega)\times \widetilde L^2(\Omega)$ (see Lemma~\ref{F:lemma: u1 is attained strongly}). Furthermore, Proposition \ref{F:prop: energy identity weakly viscous nonlocal wave} ensures that $u$ satisfies the asserted energy identity. Therefore, to complete the proof, it remains to show that
        \begin{equation}
        \label{F:eq: continuity assertion}
            u\in C([0,T];\widetilde{H}^s(\Omega))\cap C^1([0,T];\widetilde{L}^2(\Omega));
        \end{equation}
        see Definition~\ref{F:def: def of solutions}. First, recall the energy identity
        \begin{equation}
        \label{F:eq: energy identity finite energy sol}
             \begin{split}
             &\|\partial_t u(t)\|_{L^2(\Omega)}^2+\|(-\Delta)^{s/2}u(t)\|_{L^2(\R^n)}^2+2 \int_0^t \|
            (-\Delta)^{r/2} \partial_t u\|_{L^2(\R^n)}^2\,d\tau\\
             &=\|u_1\|_{L^2(\Omega)}^2+\|(-\Delta)^{s/2}u_0\|_{L^2(\R^n)}^2 -2 \re \int_0^t \langle qu,\partial_t u\rangle_{L^2(\Omega)}\,d\tau\\
            &\quad +2\re \int_0^t\langle f,\partial_t u\rangle\,d\tau
        \end{split}
        \end{equation}
        for all $0\leq t\leq T$. By absolute continuity of the integrals, this shows that
        \[
            \zeta(t)\vcentcolon = \|\partial_t u(t)\|_{L^2(\Omega)}^2+\|(-\Delta)^{s/2}u(t)\|_{L^2(\R^n)}^2\geq 0
        \]
        is a continuous function. In addition, let us define
        \[
            \Delta\zeta(t,\tau)\vcentcolon = \|\partial_t u(t)-\partial_t u(\tau)\|_{L^2(\Omega)}^2+\|(-\Delta)^{s/2}u(t)-(-\Delta)^{s/2}u(\tau)\|_{L^2(\R^n)}^2
        \]
        for all $0\leq t,\tau\leq T$. We have
        \[
        \begin{split}
            \Delta\zeta(t,\tau)&=\zeta(t)+\zeta(\tau)-2\re\langle \partial_t u(t),\partial_t u(\tau)\rangle_{L^2(\Omega)}\\
            &\quad -2\re\langle (-\Delta)^{s/2}u(t),(-\Delta)^{s/2}u(\tau)\rangle_{L^2(\R^n)}.
        \end{split}
        \]
        Using the continuity of $\zeta$ and
        \[
            (u,\partial_t u)\in C_w([0,T];\widetilde{H}^s(\Omega))\times C_w([0,T];\widetilde{L}^2(\Omega)),
        \]
        we deduce
        \[
            \Delta\zeta(t,\tau)\to 2\zeta(t)-2 \|\partial_t u(t)\|_{L^2(\Omega)}^2-2\|(-\Delta)^{s/2}u(t)\|_{L^2(\R^n)}^2=0
        \]
        as $\tau\to t$. This completes the proof.
    \end{proof}

    \subsection{Uniqueness of weak finite energy solutions}
    \label{F:sec: uniqueness finite energy sol}

    In this section, we prove uniqueness of weak finite energy solutions to \eqref{F:eq: existence weakly viscous case}. The applied method goes back to O.~A.~Ladyzhenskaya \cite{Ladyzhenskaya}. 

    \begin{theorem}
    \label{F:prop: uniqueness}
        Let $\Omega\subset\R^n$ be a bounded Lipschitz domain, $T>0$, and $0\leq r<s<\infty$. Suppose that $q\in\mathcal{Q}_s$ has integrability exponent $2\leq p\leq \infty$, and
        \[
            (u_0,u_1,f)\in \widetilde{H}^s(\Omega)\times \widetilde{L}^2(\Omega)\times L^2(0,T;H^{-r}(\Omega)).
        \]
       Then problem \eqref{F:eq: existence weakly viscous case} has at most one weak finite energy solution.
    \end{theorem}

    \begin{proof}
        Since \eqref{F:eq: existence weakly viscous case} is a linear PDE, it is enough to demonstrate that if $u$ is a weak finite energy solution of
        \begin{equation}
        \label{F:eq: PDE for uniqueness}
            \begin{cases}
				( \partial_t^2 +(-\Delta)^r \partial_t +(-\Delta)^s+q) u=0 &\text{ in }\Omega_T,\\
				u=0     &\text{ in }(\Omega_e)_T,\\
				u(0)=0,\, \partial_t u(0)=0  &\text{ in }\Omega,
			\end{cases}
        \end{equation}
        then $u=0$. Let us fix some $0<t_0<T$ and introduce the function
       \begin{equation}
       \label{F:eq: test function}
           \varphi(t)=\begin{cases}
               -\int_{t}^{t_0}u(\tau)\,d\tau,&\text{ for }0\leq t\leq t_0,\\
               0,&\text{ for } t_0\leq t\leq T.
           \end{cases}
       \end{equation}
        Let us note that $\varphi$ satisfies
      \begin{equation}
      \label{F:eq: regularity test function uniqueness}
        \varphi\in H^1(0,T;\widetilde{H}^s(\Omega))\cap C([0,T];\widetilde{H}^s(\Omega))\text{ with }
      \partial_t\varphi=u\chi_{[0,t_0]}.
      \end{equation}
      The condition $\varphi\in C([0,T];\widetilde{H}^s(\Omega))$ is easily seen by fundamental theorem of calculus. Let us shortly give the argument for $\partial_t \varphi \in L^2(0,T;\widetilde H^s(\Omega))$. Since $\varphi\in H^1(0,t_0;\widetilde H^s(\Omega))$, by usual density results, we know that there exists $(\varphi_n)_{n\in\N}\subset C^{\infty}([0,t_0];\widetilde{H}^s(\Omega))$ such that $\varphi_n\to\varphi$ in $H^1(0,t_0;\widetilde{H}^s(\Omega))$ as $n\to\infty$. Thus, for any $\eta\in C_c^{\infty}((0,T))$ we may compute
      \[
      \begin{split}
\int_0^T\varphi(t)\eta'(t)\,dt&=\int_0^{t_0}\varphi(t)\eta'(t)\,dt\\
&=\lim_{n\to\infty}\int_0^{t_0}\varphi_n(t)\eta'(t)\,dt\\
          &=\lim_{n\to\infty}\left(\varphi(t_0)\eta(t_0)-\int_0^{t_0}\varphi'_n(t)\eta(t)\,dt\right)\\
          &=\varphi(t_0)\eta(t_0)-\int_0^{t_0}\varphi'(t)\eta(t)\,dt\\
          &=-\int_0^{t_0}\varphi'(t)\eta(t)\,dt.
      \end{split}
      \]
      In the above computation we used the fundamental theorem of calculus and for the convergence of the boundary value the embedding 
      \[
        H^1(0,t_0;\widetilde{H}^s(\Omega))\hookrightarrow C([0,t_0];\widetilde{H}^s(\Omega)).
      \]
      Hence, this establishes $\partial_t \varphi=u\chi_{[0,t_0]}$, and the properties \eqref{F:eq: regularity test function uniqueness} follow.
      Thus, from \cite[Appendix B]{LZ2023unique} we may infer that $\varphi$ can be approximated by functions in $C_c^{\infty}(\Omega\times [0,T))$ with respect to the norm of the space $H^1(0,T;\widetilde{H}^s(\Omega))$. So, from \eqref{F:eq: weak formualation} and \eqref{F:eq: regularity test function uniqueness} we deduce that
       \begin{equation}
       \label{F:eq: energy identity for uniqueness}
       \begin{split}
                0=&\re\int_0^{t_0}\left(-\langle \partial_t u,u\rangle_{L^2(\Omega)} +\langle (-\Delta)^{r/2} \partial_t u,(-\Delta)^{r/2}\varphi\rangle_{L^2(\R^n)}\right)\,dt\\
            &+\re\int_0^{t_0}\left(\langle (-\Delta)^{s/2}\partial_t\varphi,(-\Delta)^{s/2}\varphi\rangle_{L^2(\R^n)}+\langle q u,\varphi\rangle_{L^2(\Omega)}\right)\,dt\\
            =&\frac{1}{2}\int_0^{t_0}\left(-\partial_t\|u\|^2_{L^2(\Omega)} +2\re\langle (-\Delta)^{r/2} \partial_t u,(-\Delta)^{r/2}\varphi\rangle_{L^2(\R^n)}\right)\,dt\\
            &+\frac{1}{2}\int_0^{t_0}\left(\partial_t\langle (-\Delta)^{s/2}\varphi,(-\Delta)^{s/2}\varphi\rangle_{L^2(\R^n)}+2\re\langle q u,\varphi\rangle_{L^2(\Omega)}\right)\,dt.
            \end{split}
       \end{equation}
       Using $u\in H^1(0,T;\widetilde{H}^r(\Omega))$, we know that $(-\Delta)^{r/2}u\in H^1(0,T;L^2(\R^n))$. So, $\varphi\in H^1(0,T;\widetilde{H}^s(\Omega))$, $\varphi(t_0)=0$ in $\widetilde{H}^s(\Omega)$, $u(0)=0$ in $\widetilde{H}^r(\Omega)$, and an integration by parts ensure
       \begin{equation}
       \label{F:eq: integration by parts damping}
       \begin{split}
        &\int_0^{t_0}\langle (-\Delta)^{r/2} \partial_t u,(-\Delta)^{r/2}\varphi\rangle_{L^2(\R^n)}\,dt\\
        &=\left.\langle (-\Delta)^{r/2} u(t),(-\Delta)^{r/2}\varphi(t)\rangle_{L^2(\R^n)}\right|_{t=0}^{t=t_0}-\int_0^{t_0}\langle (-\Delta)^{r/2} u,(-\Delta)^{r/2} \partial_t\varphi\rangle_{L^2(\R^n)}\,dt\\
        &=-\int_0^{t_0}\|(-\Delta)^{r/2} u\|^2_{L^2(\R^n)}\,dt.
        \end{split}
       \end{equation}
       Thus, we arrive at the identity
       \[
       \begin{split}
           &\|u(t_0)\|_{L^2(\Omega)}^2+\|(-\Delta)^{s/2}\varphi(0)\|_{L^2(\R^n)}^2+2\int_0^{t_0}\|(-\Delta)^{r/2}u\|_{L^2(\R^n)}^2\,dt\\
           &=2\re\int_0^{t_0}\langle q u,\varphi\rangle_{L^2(\Omega)}\, dt.
       \end{split}
       \]
       Thus, Lemma~\ref{F:lemma: cond on q} ensures that
       \[
       \begin{split}
         &\|u(t_0)\|_{L^2(\Omega)}^2+\|(-\Delta)^{s/2}\varphi(0)\|_{L^2(\R^n)}^2+2\int_0^{t_0}\|(-\Delta)^{r/2}u\|_{L^2(\R^n)}^2\,dt\\
           &\leq C\int_0^{t_0}\|q(t)\|_{L^p(\Omega)}\|u(t)\|_{L^2(\Omega)}\|(-\Delta)^{s/2}\varphi(t)\|_{L^2(\R^n)}\,dt\\
           &\leq C\int_0^{t_0}\|q(t)\|_{L^p(\Omega)}(\|u(t)\|_{L^2(\Omega)}^2+\|(-\Delta)^{s/2}\varphi(t)\|_{L^2(\R^n)}^2)\,dt.
       \end{split}
       \]
       Next, let us define
       \[
        \psi(t)\vcentcolon = \int_0^t u(\tau)\,d\tau.
       \]
       Then, the previous estimate becomes
       \begin{equation}
       \begin{split}
           &\|u(t_0)\|_{L^2(\Omega)}^2+\|(-\Delta)^{s/2}\psi(t_0)\|_{L^2(\R^n)}^2+2\int_0^{t_0}\|(-\Delta)^{r/2}u\|_{L^2(\R^n)}^2\,dt\\
           &\leq C\int_0^{t_0}\|q(t)\|_{L^p(\Omega)}(\|u(t)\|_{L^2(\Omega)}^2+\|(-\Delta)^{s/2}(\psi(t_0)-\psi(t))\|_{L^2(\R^n)}^2)\,dt\\
           &\leq C_0\|q\|_{L^1(0,t_0;L^p(\Omega))}\|(-\Delta)^{s/2}\psi(t_0)\|_{L^2(\R^n)}^2\\
           &\quad+ C\int_0^{t_0}\|q(t)\|_{L^p(\Omega)}(\|u(t)\|_{L^2(\Omega)}^2+\|(-\Delta)^{s/2}\psi(t)\|_{L^2(\R^n)}^2)\,dt.
       \end{split}
       \end{equation}
       If $0<T_0\leq T$ is chosen such that
       \[
            \|q\|_{L^1(0,T_0;L^p(\Omega))}\leq 1/(2C_0),
       \]
       then we get
        \begin{equation}
            \begin{split}
                &\|u(t_0)\|_{L^2(\Omega)}^2+\|(-\Delta)^{s/2}\psi(t_0)\|_{L^2(\R^n)}^2+2\int_0^{t_0}\|(-\Delta)^{r/2}u\|_{L^2(\R^n)}^2\,dt\\
                &\leq C\int_0^{t_0}\|q(t)\|_{L^p(\Omega)}(\|u(t)\|_{L^2(\Omega)}^2+\|(-\Delta)^{s/2}\psi(t)\|_{L^2(\R^n)}^2)\,dt
            \end{split}
        \end{equation}
       for all $0\leq t_0\leq T_0$ and some constant $C>0$ depending only on $p$, $n$, $s$, and $\Omega$. Thus, Gronwall's inequality implies $u=0$ on $[0,T_0]$. If we can take $T_0=T$ this finishes the proof, which is in particular true if $\|q\|_{L^1(0,T;L^p(\Omega))}$ is sufficiently small. On the other hand, if $T_0<T$, then we observe that $u$ satisfies 
       \[
        \begin{cases}
				( \partial_t^2 +(-\Delta)^r \partial_t +(-\Delta)^s+q) u=0 &\text{ in }\Omega\times (T_0,T),\\
				u=0     &\text{ in }\Omega_e\times (T_0,T),\\
				u(T_0)=0,\, \partial_t u(T_0)=0  &\text{ in }\Omega,
			\end{cases}
       \]
       and we may repeat the above argument to find $0<T_1\leq T$ such that $u=0$ on $[T_0,T_1]$. Thus, we find an increasing sequence $(T_n)_{n\in \Z_{\geq -1}}\subset [0,T]$ with $T_{-1}=0$, $T_n\to T$ as $n\to\infty$ and $u=0$ on $[T_n,T_{n+1}]$ for all $n\in\Z_{\geq -1}$. Hence, we may conclude that $u=0$ on $[0,T]$. In fact, as $q\in L^{\infty}(0,T;L^p(\Omega))$ and $C_0>0$ only depends on $n,s,p$, and $\Omega$, it is sufficient to perform this iteration finitely many times to conclude the desired assertion.
    \end{proof}

    \subsection{Nonhomogeneous linear weakly viscous equations}
    \label{F:sec: nonhomgeneous weakly viscous}

    We use the following definition of solutions for weakly viscous equations with non-zero exterior conditions.

    \begin{definition}
    \label{F:def: def of weak sol nonhomogeneous}
        Let $\Omega\subset\R^n$ be a bounded Lipschitz domain, $T>0$, and $0\leq r<s<\infty$. Furthermore, assume that $q\in\mathcal{Q}_s$ with integrability exponent $2\leq p\leq \infty$. Suppose that
        \[
            (u_0,u_1,\varphi)\in H^s(\R^n)\times L^2(\R^n)\times C^2([0,T];H^{2s}(\R^n))
        \]
        satisfies the compatibility condition
        \begin{equation}
        \label{F:eq: compatibility cond}
            \varphi(0)=u_0\text{ and }\partial_t\varphi(0)=u_1\text{ a.e. in }\Omega_e.
        \end{equation}
        We say that a function 
        \[
            u\in L^2(0,T;H^s(\R^n))\cap H^1(0,T;H^r(\R^n))
        \]
        is a \emph{weak finite energy solution} of 
        \begin{equation}
			\label{F:eq: existence weakly viscous case nonhom}
			\begin{cases}
				( \partial_t^2 +(-\Delta)^r \partial_t +(-\Delta)^s+q) u=0 &\text{ in }\Omega_T,\\
				u=\varphi     &\text{ in }(\Omega_e)_T,\\
				u(0)=u_0,\, \partial_t u(0)=u_1  &\text{ in }\Omega,
			\end{cases}
	\end{equation}
        if
        \begin{enumerate}[(a)]
            \item  $u$ solves \eqref{F:eq: existence weakly viscous case nonhom} in the weak sense; see equation~\eqref{F:eq: calc for distr sol}, 
            \item $u=\varphi$ in $(\Omega_e)_T$,
            \item and $u(0)=u_0$ in $H^r(\R^n)$.
        \end{enumerate}
        If $u$ belongs to the space $C([0,T];H^s(\R^n))\cap C^1([0,T];L^2(\R^n))$, and $(u,\partial_t u)(0)=(u_0,u_1)$ in $H^s(\R^n)\times L^2(\R^n)$, then we say $u$ is a \emph{finite energy solution}.
    \end{definition}

    From the results established in the previous sections, we deduce the following well-posedness result.

    \begin{theorem}[Well-posedness non-homogeneous, weakly viscous wave equations]
     \label{F:thm: unique solvability of weakly viscous wave eq nonhom}
         Let $\Omega\subset\R^n$ be a bounded Lipschitz domain, $T>0$, and $0\leq r<s<\infty$. Furthermore, assume that $q\in\mathcal{Q}_s$ with integrability exponent $2\leq p\leq \infty$. Suppose that
         \[
            (u_0,u_1,\varphi)\in H^s(\R^n)\times L^2(\R^n)\times C^2([0,T];H^{2s}(\R^n))
         \]
        satisfies the compatibility condition \eqref{F:eq: compatibility cond}. 
        
        Then problem \eqref{F:eq: existence weakly viscous case nonhom} has a unique finite energy solution $u$ such that $u-\varphi$ satisfies the energy identity \eqref{F:eq: energy identity weakly viscous nonlocal wave eq continuous potential}.
     \end{theorem}

     \begin{proof}
         Let us observe that $u$ is a finite energy solution of \eqref{F:eq: existence weakly viscous case nonhom} if and only if $v=u-\varphi$ is a finite energy solution of \eqref{F:eq: existence weakly viscous case} with $(u_0,u_1,f)$ replaced by
         \begin{equation}
         \label{F:eq: trafo with exterior cond}
            \begin{split}
                &u_0-\varphi(0)\in\widetilde{H}^s(\Omega),\,u_1-\partial_t \varphi(0)\in \widetilde{L}^2(\Omega),\\
                &f=-\chi_{\Omega} (\partial_t^2\varphi+(-\Delta)^r \partial_t \varphi+(-\Delta)^s \varphi+q\varphi)\in L^2(0,T;\widetilde{L}^2(\Omega)).
            \end{split}
         \end{equation}
        Here we use that, for every $w\in H^s(\R^n)$, one has $q\,w|_\Omega \in L^2(\Omega)$.
        Indeed, this follows from the proof of \cite[Theorem~3.1]{zimmermann2024-viscous-wave}.
         Thus, the assertions are direct consequences of Theorem~\ref{F:prop: existence of weak solutions} and Theorem~\ref{F:prop: uniqueness}.
     \end{proof}
     
\section{Unique Continuation and Entanglement for fractional Laplacians}
\label{R:sec: UCP and entanglement}

In this section we recall the unique continuation principle for the fractional Laplacian and the recently established entanglement principle.

\subsection{Unique continuation}
\label{R:subsec: UCP}

A fundamental ingredient in the proof of the Runge approximation theorem is the unique continuation property (UCP) for fractional Laplacians.

\begin{proposition}[Unique continuation for fractional Laplacians]
\label{R:prop:UCP}
Let $s>0$ be non-integer, $1<p<\infty$, and $r\in\mathbb{R}$. If $u\in H^{r,p}(\R^n)$ satisfies
\[
u = (-\Delta)^s u = 0
\quad\text{in a nonempty open set }V\subset\R^n,
\]
then $u\equiv0$ in $\R^n$.
\end{proposition}

For $0<s<1$, $p=2$, and $r\in\R$, Proposition~\ref{R:prop:UCP} was first established in \cite[Theorem~1.2]{GSU20}. In this range of fractional powers $s$, the fractional Laplacian admits the singular integral representation
\[
(-\Delta)^s u(x)
= C_{n,s}\,\mathrm{p.v.}\!\int_{\mathbb{R}^n}
\frac{u(x)-u(y)}{|x-y|^{n+2s}}\,dy,
\]
for sufficiently regular functions $u$. For higher non-integer orders $s>1$, one may write
\[
s=m+\sigma,
\qquad
m\in\mathbb N,\quad 0<\sigma<1,
\]
and combine the UCP for $(-\Delta)^\sigma$ with the classical unique continuation property for integer powers of the Laplacian. The UCP was extended to the full scale of Bessel potential spaces $H^{r,p}(\R^n)$ in \cite[Theorem~2.2]{KRZ-2023}. Eventually, it has been shown in \cite[Theorem~3.1]{uhlmann2026fractionalcalderonproblemlinftyconductivities} that the UCP property remains valid for functions $u\in L^1_{\rm loc}(\R^n)$ satisfying 
\[
    \int_{\R^n}\frac{|u(x)|}{1+|x|^{n+2s}}\,dx<\infty.
\]

\subsection{Entanglement principle}
\label{R:subsec: entanglement}

The following result allows one to decouple linear combinations of fractional Laplacians of different orders.

  \begin{proposition}[Entanglement principle in Bessel potential spaces]
\label{R:prop: entanglement Bessel}
Assume $n\geq2$.
Let $1<p<\infty$, $t\in\R$, and let
\[
u_j\in H^{t,p}(\R^n),
\qquad j=1,\dots,N,
\]
satisfy the super-exponential decay condition
\[
|\langle u_j,\phi\rangle|
\le
Ce^{-\rho R^\gamma}
\|\phi\|_{H^{-t,p'}(\R^n)}
\]
for all $\phi\in C_c^\infty(\R^n \setminus \overline B_R)$, where $C,\rho>0$, $\gamma>1$, and $p'=p/(p-1)$.

Assume that
\[
u_j=0 \quad \text{in } O,
\qquad j=1,\dots,N,
\]
for some nonempty open set $O\subset\R^n$, and
\[
\sum_{j=1}^N b_j(-\Delta)^{s_j}u_j=0
\quad\text{in }O,
\]
where the coefficients $(b_j)_{j=1,\ldots,N}\subset \C$ are nonzero and the exponents $(s_j)_{j=1,\ldots,N}\subset\R_+\setminus\mathbb N$ are monotonically increasing and satisfy 
    \begin{equation}
\label{R:eq: cond fractional exponents}
\begin{cases}
s_k-s_j\notin\mathbb N,
& \text{for all }k>j,
\quad \text{if } n\in2\mathbb N,
\\[1mm]
s_k-s_j\notin\frac12\mathbb N,
& \text{for all }k>j,
\quad \text{if } n\in 2\mathbb N+1.
\end{cases}
\end{equation}
Then $u_j\equiv 0$ in $\R^n$ for all $j=1,\ldots, N$.
\end{proposition}

We remark that this result has been established in \cite[Theorem~1.2]{feizmohammadi2024entanglementprinciplefractionallaplacian} for $p=2$, which in turn builds up on ideas from \cite{feizmohammadi2024calder}. The extension to the full scale of Bessel potential spaces $H^{t,p}(\R^n)$, $1<p<\infty$, follows by adapting the arguments of \cite{feizmohammadi2024entanglementprinciplefractionallaplacian}. For the convenience of the reader, we include a proof in Appendix~\ref{R:sec: proof entanglement principle}.

\section{Weakly viscous equations with constant coefficients}
\label{R:sec: Runge for constant coeff}

In this section we establish several Runge approximation properties for weakly viscous nonlocal wave equations in different function spaces. Throughout this section, $\Omega\subset\R^n$ denotes a bounded Lipschitz domain, $W\subset\Omega_e$ is a fixed nonempty open measurement set, $T>0$ is a finite time-horizon, $0\leq r<s$, and $q\in \mathcal{Q}_s^p$.

We define the associated \emph{Runge set} by
\begin{equation}
\label{R:eq: Runge set}
\mathscr R_W
\vcentcolon=
\{u_\varphi-\varphi\,;\,\varphi\in C_c^\infty(W_T)\},
\end{equation}
where $u_\varphi$ denotes the unique finite energy solution of
\begin{equation}
\label{R:eq: Wave eq for Runge}
\begin{cases}
(\partial_t^2+(-\Delta)^r\partial_t+(-\Delta)^s+q)u=0
&\text{in }\Omega_T,\\
u=\varphi
&\text{in }(\Omega_e)_T,\\
u(0)=0,\ \partial_tu(0)=0
&\text{in }\Omega.
\end{cases}
\end{equation}
We refer to Remark~\ref{F:thm: unique solvability of weakly viscous wave eq nonhom} for the well-posedness of nonhomogeneous exterior value problems.

We prove density in $L^2(0,T;\widetilde H^r(\Omega))$ without smallness of $q$ and deduce Lebesgue-space density below. Stronger spatial and space--time approximation is treated in \cite{zimmermann-paperII}.

 \subsection{Runge approximation in \texorpdfstring{$L^2(0,T;\widetilde{H}^r(\Omega))$}{space sobolev}}
\label{R:sec: runge in Hr}

\begin{theorem}[Runge approximation in $L^2(0,T;\widetilde{H}^r(\Omega))$]
\label{R:thm: Runge approx in Hr}
Assume that the exponents $r,s$ satisfy one of the following conditions:
\begin{enumerate}[(i)]
\item\label{R:cond 1 exp Runge}
$s\in\N$, $r\in\R_+\setminus\N$;

\item\label{R:cond 2 exp Runge}
$s\in\R_+\setminus\N$, $r\in\N_0$;

\item\label{R:cond 3 exp Runge}
Assume $n\geq2$.
$r,s\in\R_+\setminus\N$ and
\begin{equation}
\begin{cases}
s-r\notin\mathbb N,
&
\quad \text{if } n\in2\mathbb N,
\\[1mm]
s-r\notin\frac12\mathbb N,
& 
\quad \text{if } n\in 2\mathbb N+1.
\end{cases}
\end{equation}
\end{enumerate}
Then the Runge set $\mathscr R_W$, defined in \eqref{R:eq: Runge set}, is dense in $L^2(0,T;\widetilde H^r(\Omega))$.
\end{theorem}

\begin{proof}
Since the Runge set $\mathscr{R}_W$ is a subspace of $L^2(0,T;\widetilde H^r(\Omega))$, by the Hahn--Banach theorem, it suffices to show that if
\[
f\in L^2(0,T;H^{-r}(\Omega))
\]
satisfies
\begin{equation}
\label{R:eq: cond Runge 1}
\langle f,u_\varphi-\varphi\rangle =0
\qquad
\text{for all }
\varphi\in C_c^\infty(W_T),
\end{equation}
then necessarily $f=0$.

Let $v$ denote the unique finite energy solution of the adjoint problem
\begin{equation}
\label{R:eq: adjoint problem first runge}
\begin{cases}
(\partial_t^2-(-\Delta)^r\partial_t+(-\Delta)^s+q)v=f
&\text{in }\Omega_T,\\
v=0
&\text{in }(\Omega_e)_T,\\
v(T)=0,\quad \partial_t v(T)=0
&\text{in }\Omega.
\end{cases}
\end{equation}
The well-posedness of \eqref{R:eq: adjoint problem first runge} is an immediate consequence of Theorem~\ref{F:prop: existence of weak solutions} after time-reversal $t\mapsto T-t$; the reversed potential is $q(x,T-t)$, which belongs to the same potential class.

For $\varphi\in C_c^\infty(W_T)$, define
\[
w\vcentcolon=u_\varphi-\varphi\quad \text{and}\quad g_\varphi=-((-\Delta)^r\partial_t+(-\Delta)^s)\varphi.
\]
Then $w$ solves
\begin{equation}
\label{R:eq: solution w}
\begin{cases}
(\partial_t^2+(-\Delta)^r\partial_t+(-\Delta)^s+q)w
= g_\varphi
&\text{in }\Omega_T,\\
w=0
&\text{in }(\Omega_e)_T,\\
w(0)=0,\quad \partial_t w(0)=0
&\text{in }\Omega.
\end{cases}
\end{equation}

Let
\[
w_\eps,\, v_\eps \in H^1(0,T;\widetilde{H}^s(\Omega))\cap H^2(0,T;H^{-s}(\Omega))
\]
denote the corresponding parabolic regularizations, that is, the unique weak solutions of the parabolically regularized problems
\begin{equation}
\label{R:eq: regularized solution w}
\begin{cases}
(\partial_t^2+\eps (-\Delta)^s\partial_t +(-\Delta)^r\partial_t+(-\Delta)^s+q)w
= g_\varphi
&\text{in }\Omega_T,\\
w=0
&\text{in }(\Omega_e)_T,\\
w(0)=0,\quad \partial_t w(0)=0
&\text{in }\Omega
\end{cases}
\end{equation}
and 
\begin{equation}
\label{R:eq: regularized adjoint problem first runge}
\begin{cases}
(\partial_t^2-\eps(-\Delta)^s\partial_t-(-\Delta)^r\partial_t+(-\Delta)^s+q)v=f
&\text{in }\Omega_T,\\
v=0
&\text{in }(\Omega_e)_T,\\
v(T)=0,\quad \partial_t v(T)=0
&\text{in }\Omega.
\end{cases}
\end{equation}

Using the regularized equations, we get
\[
\begin{aligned}
0
&=
\langle f,w\rangle
\\
&=
\lim_{\eps\to0}
\langle f,w_\eps\rangle
\\
&=
\lim_{\eps\to0}
\Big\langle
(\partial_t^2
-\eps(-\Delta)^s\partial_t
-(-\Delta)^r\partial_t
+(-\Delta)^s+q)v_\eps,
w_\eps
\Big\rangle
\\
&=
\lim_{\eps\to0}
\Big\langle
v_\eps,
(\partial_t^2
+\eps(-\Delta)^s\partial_t
+(-\Delta)^r\partial_t
+(-\Delta)^s+q)w_\eps
\Big\rangle
\\
&=
-\lim_{\eps\to0}
\Big\langle
v_\eps,
((-\Delta)^r\partial_t+(-\Delta)^s)\varphi
\Big\rangle
\\
&=
-
\Big\langle
v,
((-\Delta)^r\partial_t+(-\Delta)^s)\varphi
\Big\rangle .
\end{aligned}
\]

More precisely, we used:
\begin{itemize}
\item
condition~\eqref{R:eq: cond Runge 1};

\item
the weak convergence
\[
w_\eps\rightharpoonup w
\quad
\text{in }
L^2(0,T;\widetilde H^s(\Omega));
\]

\item
the fact that
\[
f\in L^2(0,T;H^{-r}(\Omega))
\subset
L^2(0,T;H^{-s}(\Omega));
\]

\item
the regularity
\[
w_\eps,v_\eps
\in
H^1(0,T;\widetilde{H}^s(\Omega))\cap H^2(0,T;H^{-s}(\Omega))
\]
for the integration by parts;
\item
and finally
\[
v_\eps\weak v
\quad
\text{in }
L^2(0,T;\widetilde L^2(\Omega)).
\]
\end{itemize}

Consequently,
\[
\langle
v,
((-\Delta)^r\partial_t+(-\Delta)^s)\varphi
\rangle
=0
\]
for all $\varphi\in C_c^\infty(W_T)$. The above convergence properties are recorded in Proposition~\ref{F:prop: limit function}. Hence,
\begin{equation}
\label{R:eq: UCP cond Runge 1}
((-\Delta)^r\partial_t-(-\Delta)^s)v=0
\quad
\text{in }W_T.
\end{equation}

We now distinguish the different cases.

\medskip

\noindent
\textbf{Case \ref{R:cond 1 exp Runge}.}
Suppose that
\[
s\in\N,
\qquad
r\in\R_+\setminus\N.
\]
Then \eqref{R:eq: UCP cond Runge 1} implies
\[
\partial_t v
=
(-\Delta)^r\partial_t v
=
0
\qquad
\text{in }W_T.
\]
Since $r$ is non-integer, Proposition~\ref{R:prop:UCP} yields $\partial_t v\equiv 0$ in $\R^n_T$. Using $v\in H^1(0,T;\widetilde{H}^r(\Omega))$ and the terminal condition $v(T)=0$, we conclude that $v\equiv 0$ in $\R^n_T$.
\medskip

\noindent
\textbf{Case \ref{R:cond 2 exp Runge}.}
If
\[
s\in\R_+\setminus\N,
\qquad
r\in\N_0,
\]
then \eqref{R:eq: UCP cond Runge 1} together with Proposition~\ref{R:prop:UCP} directly implies $v\equiv 0$ in $\R^n_T$.
 
\medskip

\noindent
\textbf{Case \ref{R:cond 3 exp Runge}.}
Assume now that
\[
r,s\in\R_+\setminus\N.
\]

If
\[
n\in2\N+1,
\qquad
s-r\notin \frac12\N,
\]
or
\[
n\in2\N,
\qquad
s-r\notin\N,
\]
then Proposition~\ref{R:prop: entanglement Bessel} applied to
\eqref{R:eq: UCP cond Runge 1} implies $v\equiv 0$ in $\R^n_T$.

Thus, in all cases, $v\equiv 0$ in $\R^n_T$.
Since $v$ weakly solves \eqref{R:eq: adjoint problem first runge}, we conclude that $f\equiv 0\quad \text{in }\R^n_T$.
This completes the proof.
\end{proof}

\section{Lebesgue approximation and recovery of the potential}
\label{sec:Lebesgue-Runge-inverse}

In this final section, we establish the unique recover potentials from the DN map associated to weakly viscous nonlocal wave equations. Our uniqueness theorem follows the same solution comparison used in
\cite[Section~4.2, proof of Theorem~1.1 for generic potentials]{zimmermann2024-viscous-wave}.
In that work the two fractional powers coincide, and the exterior equation
reduces, after unique continuation, to an ordinary differential equation in
time. Here the powers satisfy $r<s$, so we instead use locality when $r$ is
an integer and the entanglement principle when both powers are noninteger.
The last step uses only Lebesgue-space Runge approximation.

\subsection{Lebesgue-space Runge approximation}
\label{subsec:Lebesgue-Runge}

An immediate consequence of the above Runge approximation theorem is the follolwing corollary.

\begin{corollary}[Lebesgue Runge approximation]\label{new:L2}
Under the assumptions of Theorem~\ref{R:thm: Runge approx in Hr}, the restrictions $u_\varphi|_{\Omega_T}$ are dense in $L^2(\Omega_T)$.
\end{corollary}
\begin{proof}
The Runge set is dense in $L^2(0,T;\widetilde H^r(\Omega))$, which is continuously and densely embedded in $L^2(\Omega_T)$. Approximate an arbitrary $L^2$ target first by a smooth compactly supported function and then by a Runge sequence. Exterior controls vanish in $\Omega_T$, so $u_\varphi-\varphi=u_\varphi$ there.
\end{proof}

\subsection{The partial Dirichlet-to-Neumann map}

Let $\Omega\subset\R^n$ be a bounded Lipschitz domain, let
$\Omega_e=\R^n\setminus\overline\Omega$, and write $A_T=A\times(0,T)$.
For a real potential $q\in\mathcal Q_s$ and
$\varphi\in C_c^\infty((\Omega_e)_T)$, denote by $u_{q,\varphi}$ the
unique finite energy solution of
\begin{equation}\label{new:inverse-forward}
\begin{cases}
\bigl(\partial_t^2+(-\Delta)^r\partial_t+(-\Delta)^s+q\bigr)u=0
    &\text{in }\Omega_T,\\
u=\varphi&\text{in }(\Omega_e)_T,\\
u(0)=0,\quad\partial_tu(0)=0&\text{in }\Omega.
\end{cases}
\end{equation}
Existence and uniqueness follow from
Theorem~\ref{F:thm: unique solvability of weakly viscous wave eq nonhom}.
In particular,
\begin{equation}\label{new:inverse-regularity}
u_{q,\varphi}\in C([0,T];H^s(\R^n))
\cap C^1([0,T];L^2(\R^n))
\cap H^1(0,T;H^r(\R^n)).
\end{equation}

For clarity, all data and pairings in this subsection are real-valued.
This entails no loss for the uniqueness result with real potentials: equality
of the complex DN maps implies equality on real test functions. Define
$\Lambda_q\colon C_c^\infty((\Omega_e)_T)\to\mathscr D'((\Omega_e)_T)$ by
\begin{equation}\label{new:DN-pairing}
\begin{split}
\langle\Lambda_q\varphi,\psi\rangle
={}&\int_0^T\int_{\R^n}
(-\Delta)^{r/2}\partial_tu_{q,\varphi}\,
(-\Delta)^{r/2}\psi\,dx\,dt\\
&+\int_0^T\int_{\R^n}
(-\Delta)^{s/2}u_{q,\varphi}\,
(-\Delta)^{s/2}\psi\,dx\,dt,
\qquad \psi\in C_c^\infty((\Omega_e)_T).
\end{split}
\end{equation}
Both integrals are well defined by \eqref{new:inverse-regularity} and the
Cauchy--Schwarz inequality. The resulting functional is continuous in the
test-function topology. By the distributional definition of the fractional
Laplacian, this is precisely
\[
\Lambda_q\varphi
=\left.\bigl((-\Delta)^r\partial_t+(-\Delta)^s\bigr)
 u_{q,\varphi}\right|_{(\Omega_e)_T}.
\]
Thus equality of the partial maps on $(W_2)_T$ means equality of
\eqref{new:DN-pairing} for all $\psi\in C_c^\infty((W_2)_T)$.
The spatial integrals in the definition are over $\R^n$, as required by
nonlocality, even though the measurement tests are supported in the exterior.

\subsection{Partial-data recovery of a time-dependent potential}

\begin{theorem}[Uniqueness of potentials]
\label{new:potential}
Let $\Omega\subset\R^n$ be a bounded Lipschitz domain, $T>0$, and let
$W_1,W_2\subset\Omega_e$ be nonempty open sets. Suppose that $0\leq r<s$
and one of the following conditions holds:
\begin{enumerate}[(i)]
\item $s\in\N$ and $r\notin\N_0$;
\item $s\notin\N$ and $r\in\N_0$;
\item $r,s\notin\N_0$, $n\geq2$, and
\[
\begin{cases}
s-r\notin\N,& n\text{ even},\\
s-r\notin\tfrac12\N,& n\text{ odd}.
\end{cases}
\]
\end{enumerate}
Let $q_j\in\mathcal Q_s$, $j=1,2$, be real-valued,
and define $\Lambda_{q_j}$ by \eqref{new:DN-pairing}.
If
\begin{equation}\label{new:partial-data-equality}
\langle(\Lambda_{q_1}-\Lambda_{q_2})\varphi,\psi\rangle=0
\end{equation}
for every $\varphi\in C_c^\infty((W_1)_T)$ and
$\psi\in C_c^\infty((W_2)_T)$, then
$q_1=q_2$ a.e.~in $\Omega_T$.
\end{theorem}

\begin{proof}
\emph{Step 1: Compare solutions with the same exterior data.}
Fix $\varphi\in C_c^\infty((W_1)_T)$ and set
\[
u_j=u_{q_j,\varphi},\qquad w=u_1-u_2.
\]
The common exterior and initial conditions imply
\[
w=0\ \text{in }(\Omega_e)_T,\qquad w(0)=\partial_tw(0)=0.
\]
Furthermore, subtracting the exterior lifts gives
\begin{equation}\label{new:difference-regularity}
w\in C([0,T];\widetilde H^s(\Omega))
\cap C^1([0,T];\widetilde L^2(\Omega))
\cap H^1(0,T;\widetilde H^r(\Omega)).
\end{equation}
In particular, both $w(t)$ and $\partial_tw(t)$ vanish in $\Omega_e$
and have support in the fixed compact set $\overline\Omega$ for almost
every $t$. By the definition of the DN map,
\eqref{new:partial-data-equality} gives
\begin{equation}\label{new:exterior-difference}
(-\Delta)^r\partial_tw+(-\Delta)^s w=0
\quad\text{in }\mathscr D'((W_2)_T).
\end{equation}

\emph{Step 2: Justify the spatial identities at almost every time.}
The distribution in \eqref{new:exterior-difference} is represented by
\[
F(t)=(-\Delta)^r\partial_tw(t)+(-\Delta)^s w(t)
\in L^2(0,T;H^{-s}(\R^n)),
\]
because $(-\Delta)^r\partial_tw\in L^2(0,T;H^{-r}(\R^n))$,
$(-\Delta)^s w\in L^2(0,T;H^{-s}(\R^n))$, and $r<s$.
Choose a nonempty open ball $V$ with $\overline V\subset W_2$.
Testing \eqref{new:exterior-difference} with $\eta(t)\zeta(x)$ yields
\[
\int_0^T\eta(t)\langle F(t),\zeta\rangle\,dt=0
\quad\text{for all }\eta\in C_c^\infty(0,T),\quad
\zeta\in C_c^\infty(V).
\]
For each fixed $\zeta$ the scalar pairing therefore vanishes almost
everywhere in time. Choose a countable subset of $C_c^\infty(V)$ dense
in $\widetilde H^s(V)$, discard the union of the corresponding null sets,
and use continuity of $F(t)$ on $\widetilde H^s(V)$.
Together with \eqref{new:difference-regularity}, this gives, outside one
null set,
\begin{equation}\label{new:slice-identities}
w(t)=\partial_tw(t)=0\ \text{in }V,
\qquad (-\Delta)^r\partial_tw(t)+(-\Delta)^s w(t)=0
\ \text{in }\mathscr D'(V).
\end{equation}
This argument is made for the fixed input $\varphi$; no common null set
for all possible exterior inputs is needed.

\emph{Step 3: Recover equality of the interior solutions.}
We distinguish the three exponent regimes in the statement.

First suppose $s\in\N$ and $r\notin\N_0$.
Since $(-\Delta)^s$ is local and $w(t)=0$ in $V$,
we have $(-\Delta)^s w(t)=0$ in $V$ for almost every $t$.
Consequently, \eqref{new:slice-identities} gives
\[
\partial_tw(t)=(-\Delta)^r\partial_tw(t)=0
\quad\text{in }V.
\]
Because $r$ is noninteger and $\partial_tw(t)\in H^r(\R^n)$,
Proposition~\ref{R:prop:UCP} implies
$\partial_tw(t)=0$ in $\R^n$ for almost every $t\in(0,T)$.
By \eqref{new:difference-regularity},
$w\in H^1(0,T;\widetilde H^r(\Omega))$.
The Bochner fundamental theorem of calculus and the common initial
displacement therefore yield
\[
w(t)=w(0)+\int_0^t\partial_\tau w(\tau)\,d\tau=0
\quad\text{in }\widetilde H^r(\Omega),
\qquad 0\leq t\leq T.
\]
Thus the time-independent component left by the exterior identity is
eliminated by $w(0)=0$. This case requires neither entanglement nor an
additional dimension restriction.

Next suppose $s\notin\N$ and $r\in\N_0$.
Then $(-\Delta)^r$ is local, including the identity operator when $r=0$.
Since $\partial_tw(t)=0$ in $V$,
\eqref{new:slice-identities} reduces to
\[
w(t)=(-\Delta)^s w(t)=0\quad\text{in }V.
\]
As $s$ is noninteger and $w(t)\in H^s(\R^n)$,
Proposition~\ref{R:prop:UCP} implies $w(t)=0$ in $\R^n$.

Finally, suppose $r,s\notin\N_0$ and the third exponent condition holds.
For almost every $t$, apply
Proposition~\ref{R:prop: entanglement Bessel} to
\[
f_1=\partial_tw(t),\qquad f_2=w(t),
\qquad s_1=r,\quad s_2=s,\quad b_1=b_2=1.
\]
Both functions belong to $L^2(\R^n)$, vanish in $V$, and have compact
support in $\overline\Omega$; hence they satisfy the required
super-exponential decay condition. Indeed, for test functions supported
outside a sufficiently large ball their pairings are zero, while for the
remaining radii the estimate follows from their $L^2$ norms by increasing
the constant. The two powers are distinct and noninteger, and their
difference satisfies exactly the dimension-dependent nonresonance
condition in the theorem. Thus the entanglement principle gives
$f_1=f_2=0$ in $\R^n$.

In all three cases we have proved
\begin{equation}\label{new:same-input-solutions}
u_{q_1,\varphi}=u_{q_2,\varphi}
\quad\text{in }\R^n_T.
\end{equation}
In particular this equality holds as an equality of space--time
distributions, so their distributional time derivatives agree as well.

\emph{Step 4: Subtract the interior equations.}
Write $u=u_{q_1,\varphi}=u_{q_2,\varphi}$.
Using \eqref{new:inverse-forward} for the two potentials and
\eqref{new:same-input-solutions} gives
\begin{equation}\label{new:potential-times-solution}
(q_1-q_2)u=0\quad\text{in }\mathscr D'(\Omega_T).
\end{equation}
Since $q_j\in\mathcal Q_s$, the admissible exponents satisfy $p_j\geq2$. As $\Omega_T$ has finite measure, we have $q_j\in L^2(\Omega_T)$. Together with $u\in L^2(\Omega_T)$, H\"older's inequality therefore gives
\[
    (q_1-q_2)u\in L^1(\Omega_T).
\]
Consequently, \eqref{new:potential-times-solution} implies $(q_1-q_2)u=0$ a.e.~in $\Omega_T$. This conclusion holds for every $\varphi\in C_c^\infty((W_1)_T)$.

\emph{Step 5: Use Lebesgue Runge approximation.}
The hypotheses on $r,s$ are among those of
Theorem~\ref{R:thm: Runge approx in Hr}. Applying
Corollary~\ref{new:L2} to the fixed potential $q_1$ and control set $W_1$,
we obtain $\varphi_k\in C_c^\infty((W_1)_T)$ such that
\[
u_k\vcentcolon=u_{q_1,\varphi_k}|_{\Omega_T}\to1
\quad\text{in }L^2(\Omega_T).
\]
The constant target is admissible because $|\Omega_T|<\infty$. Setting $a=q_1-q_2$ and using
\eqref{new:potential-times-solution} for each $\varphi_k$, we find
\[
\|a\|_{L^1(\Omega_T)}
=\|a(1-u_k)\|_{L^1(\Omega_T)}
\leq\|a\|_{L^2(\Omega_T)}
      \|1-u_k\|_{L^2(\Omega_T)}\to0.
\]
Therefore $q_1=q_2$ almost everywhere in $\Omega_T$.
\end{proof}

\begin{remark}[Role of the measurement sets and the forward theory]
The set $W_2$ is used for unique continuation or entanglement, whereas $W_1$ is used for Runge approximation. The proof does not require
time-reversal invariance of the potentials, disjoint measurement sets, or a smallness assumption on $q_j$. It uses the finite-energy forward theory and Lebesgue Runge approximation, not the very weak theory needed for the stronger approximation results of
\cite{zimmermann-paperII}.
\end{remark}

\appendix
    \section{Proofs of auxiliary lemmas}
    \label{F:appendix: auxiliary lemmas}

    \subsection{Proof of Gronwall's inequality}
    \label{F:appendix: Proof of Gronwall}

    In this section, we give the proof of Gronwall's inequality as stated in Lemma~\ref{F:lemma: Gronwall}.

   \begin{proof}[Proof of Gronwall's inequality]
Set
\[
M(t)\vcentcolon=\int_0^t\mu(\tau)\varphi(\tau)\,d\tau .
\]
Then $M$ is absolutely continuous, $M(0)=0$, and
\[
M'(t)=\mu(t)\varphi(t)
\]
for a.e. $0<t<T$. By \eqref{F:eq: gronwall assumpt},
\[
M'(t)
\leq
\mu(t)\alpha(t)+\mu(t)M(t)
\]
for a.e. $0<t<T$. Hence
\[
\frac{d}{dt}
\left(
M(t)\exp\left(-\int_0^t\mu(\tau)\,d\tau\right)
\right)
\leq
\mu(t)\alpha(t)
\exp\left(-\int_0^t\mu(\tau)\,d\tau\right).
\]
Integrating from $0$ to $t$, we obtain
\[
M(t)
\leq
\int_0^t
\mu(s)\alpha(s)
\exp\left(\int_s^t\mu(\tau)\,d\tau\right)\,ds .
\]
Since $\alpha$ is nondecreasing, $\alpha(s)\leq \alpha(t)$ for $0\leq s\leq t$. Therefore,
\[
M(t)
\leq
\alpha(t)
\int_0^t
\mu(s)
\exp\left(\int_s^t\mu(\tau)\,d\tau\right)\,ds .
\]
Using
\[
\frac{d}{ds}
\exp\left(\int_s^t\mu(\tau)\,d\tau\right)
=
-\mu(s)
\exp\left(\int_s^t\mu(\tau)\,d\tau\right),
\]
we get
\[
\int_0^t
\mu(s)
\exp\left(\int_s^t\mu(\tau)\,d\tau\right)\,ds
=
\exp\left(\int_0^t\mu(\tau)\,d\tau\right)-1.
\]
Thus,
\[
M(t)
\leq
\alpha(t)
\left[
\exp\left(\int_0^t\mu(\tau)\,d\tau\right)-1
\right].
\]
Inserting this into \eqref{F:eq: gronwall assumpt} gives
\[
\varphi(t)
\leq
\alpha(t)+M(t)
\leq
\alpha(t)
\exp\left(\int_0^t\mu(\tau)\,d\tau\right)
\]
for a.e. $0<t<T$, as claimed.
\end{proof}

    \subsection{Generalized vector-valued Sobolev}
    \label{F:appendix: generlaized vector-valued Sobolev}

    In this section, we introduce generalized vector-valued Sobolev spaces, which we need for the proof of the energy identity (Proposition~\ref{F:prop: energy identity weakly viscous nonlocal wave}).

      Let $X,Y$ be Hilbert spaces such that $X\hookrightarrow Y$ and the inclusion is dense. For $m\in\N$ and $-\infty\leq a\leq b\leq\infty$, we define
        \begin{equation}
        \label{F:eq: generalized Sobolev}
            W^{m}(a,b;X,Y)\vcentcolon = \{u\in L^2(a,b;X)\,;\,\partial_t^m u\in L^2(a,b;Y)\}.
        \end{equation}
        Under the norm
        \[
            \|u\|_{W^m(a,b;X,Y)}\vcentcolon = (\|u\|_{L^2(a,b;X)}^2+\|\partial_t^m u\|_{L^2(a,b;Y)}^2)^{1/2},
        \]
        the space $W^{m}(a,b;X,Y)$ becomes a Hilbert space. Similarly to the spaces $H^k(a,b;X)$, we have the following density result.
        \begin{claim}
        \label{F:claim: approximation higher order}
            The space $\test([a,b];X)$ is dense in $W^{m}(a,b;X,Y)$.
        \end{claim}

        \begin{proof}[Proof of Claim~\ref{F:claim: approximation higher order}]
        Without loss of generality, we may assume that $(a,b)=(0,T)$ for some $0<T<\infty$. Let $u\in W^{m}(0,T;X,Y)$ and suppose $\eta_0,\eta_T\in C_c^{\infty}(\R)$ are cutoff functions such that $\eta_0+\eta_T=1$ on $[0,T]$, $\eta_0$ vanishes in a neighborhood of $t=T$ and $\eta_T$ vanishes in a neighborhood of $t=0$. Then, we can write
        \[
            u=u_0+u_T\text{ with }\begin{cases}
                u_0=\eta_0 u\in W^m(0,\infty;X,Y),\\
                u_T=\eta_T u\in W^m(-\infty,T;X,Y).
            \end{cases}
        \]
        So, the assertion reduces to showing $v\in W^m(0,\infty;X,Y)$ and $w\in W^m(-\infty,T;X,Y)$ can be approximated by functions in $\test([0,\infty);X)$ and $\test((-\infty,T];X)$, respectively, but this is known by \cite[p.~496]{DautrayLionsVol5} and so the claim follows.
    \end{proof}
        
        Next, we discuss the higher-order extension of functions in $W^m(0,T;X,Y)$ to functions in $W^m(-\infty,\infty;X,Y)$. To this end, assume that the sequence of real numbers $(\lambda_j)_{j=1,\ldots,m}$ satisfies the conditions
        \begin{equation}
        \label{F:eq: def of coeff}
            \sum_{j=1}^{m}(-j)^k \lambda_j=1
        \end{equation}
        for each $k=0,1,\ldots,m-1$. Note that a unique solution $\lambda=(\lambda_1,\ldots,\lambda_{m})$ exists, since the related matrix $\Lambda$ is of Vandermond type with nonzero determinant. More concretely, it is given by
        \begin{equation}
            \Lambda =\begin{pmatrix}
                1 & 1 &\cdots & 1 & 1\\
                -1 & -2 &\cdots & -(m-1)& -m\\
                 &  &\cdots & & \\
                (-1)^{m-1} & (-2)^{m-1} &\cdots & (-(m-1))^{m-1} & (-m)^{m-1}
            \end{pmatrix}
        \end{equation}
        and its determinant is
        \[
            \det\Lambda =\prod_{0\leq i<j\leq m-1}(x_j-x_i)=\prod_{0\leq i<j\leq m-1}(-1)(j-i)\neq 0,
        \]
        where we set $x_j=-(j+1)$. 
        
        For any $u\in \test([0,T];X)$ we may now define the higher-order reflection at $t=0$ as the map $P^m_0(u)\colon [-T_m,T]\to X$ satisfying
        \begin{equation}
        \label{F:eq: higher order reflection at t=0}
             P^m_0(u) (t)\vcentcolon =
            \begin{cases}
                u(t)&\text{ for }0\leq t\leq T,\\
                \sum_{j=1}^{m}\lambda_j u(-jt)&\text{ for }-T_m\leq t\leq 0,
            \end{cases}
        \end{equation}
        where
        \[
            T_m=\begin{cases}
                T\text{ for }m=1,\\
                T/(m-1)\text{ for }m\geq 2.
            \end{cases}
        \]
        We observe that the coefficients $\lambda_j$ are chosen precisely in such a way that 
        \begin{equation}
        \label{F:eq: equal k der}
            \left.\partial_t^k\right|_{t=0}\sum_{j=1}^{m}\lambda_j u(-jt)=\sum_{j=1}^{m}(-j)^k \lambda_j \partial_t^k u(-jt)\Big|_{t=0}=\partial_t^k u(0)
        \end{equation}
        for all $0\leq k\leq m-1$. Therefore, we deduce that $P^m_0(u)\in C^{m-1}([-T_m,T];X)$ and we may compute
        \begin{equation}
        \label{F:eq: formula for mth derivative}
            \partial_t^k P^m_0(u)(t)= 
            \begin{cases}
            \partial_t^k u(t)&\text{ for }0\leq t\leq T,\\
                \sum_{j=1}^{m}(-j)^k\lambda_j \partial_t^k u(-jt)&\text{ for }-T_m\leq t< 0
            \end{cases}
        \end{equation}
        for all $0\leq k\leq m$. Hence, it is easily seen that $P^m_0$ induces a continuous map
        \begin{equation}
            P^m_0\colon \test([0,T];X)\to W^m(-T_m,T;X,Y)
        \end{equation}
        when $\test([0,T];X)$ is endowed with the topology of $W^m(0,T;X,Y)$. Thus, by Claim  \ref{F:claim: approximation higher order} it can be extended to a continuous map
        \begin{equation}
            P^m_0\colon W^m(0,T;X,Y)\to W^m(-T_m,T;X,Y).
        \end{equation}
        Moreover, it has the property that $P^m_0(u)|_{(0,T)}=u$ and, thus, is an extension operator from $(0,T)$ to $(-T_m,T)$.
        
        Similarly to the reflection at $t=0$, we can introduce the higher-order reflection $P^m_T$ at $t=T$ with the same properties. Eventually, by concatenation of the extension operators $P^m_0$ and $P^m_T$ together with an additional truncation, we get a continuous extension operator
        \begin{equation}
        \label{F:eq: extension operator}
            P^m\colon W^m(0,T;X,Y)\to W^m(-\infty,\infty;X,Y).
        \end{equation}
        Moreover, the support of the extension $P^m(u)$ is contained in $(-T_m,T+T_m)$.
        
        Clearly, the above construction also works for the spaces $W^m(a,b;X,Y)$ with $-\infty\leq a\leq b\leq\infty$. Hence, we have proved the following result.

        \begin{lemma}[Extension operator]
        \label{F:lemma: extension operator}
            Let $m\in\N$, $-\infty\leq a\leq b\leq\infty$, and suppose that $X,Y$ are two Hilbert spaces such that the inclusion $X\hookrightarrow Y$ is dense. Then, for any $\Delta>0$, there exists a continuous extension operator 
            \[
                P^m\colon W^m(a,b;X,Y)\to W^m (-\infty,\infty;X,Y)
            \]
            such that $P^m(u)|_{(a,b)}=u$ and $\supp(P^m(u))\subset [a-\Delta,b+\Delta]$.
        \end{lemma}

         Moreover, let us next specialize to the case $m=2$, $X=V$, $Y=V'$, and there exists a Hilbert space $H$ such that 
         \[
            V\hookrightarrow H\hookrightarrow V'
         \]
         forms a Gelfand triple. Then \cite[p.~480]{DautrayLionsVol5} yields that
         \begin{equation}
         \label{F:eq: interpolant 1/2}
             (X,Y)_{1/2}=H.
         \end{equation}
        Here, $(X,Y)_\theta$, $0\leq \theta\leq 1$, denotes the holomorphic interpolant of the spaces $X$ and $Y$, which satisfies
         \[
            X\hookrightarrow (X,Y)_\theta \hookrightarrow Y,
         \]
         with the first inclusion being dense, and
         \begin{equation}
         \label{F:eq: interpolation estimate}
             \|u\|_{(X,Y)_{\theta}}\lesssim \|u\|_X^{1-\theta}\|u\|_Y^{\theta}
         \end{equation}
         for all $u\in X$.

          Furthermore, \cite[Chapter XVIII, \S~1, Corollary~3]{DautrayLionsVol5}, together with \eqref{F:eq: interpolant 1/2}, ensures that for any $u\in W^2(-\infty,\infty;V,V')$ we have $\partial_t u\in L^2(\R;H)$ satisfying
          \begin{equation}
          \label{F:eq: estimate for first derivative}
              \|\partial_t u\|_{L^2(\R;H)}\leq \|u\|_{L^2(\R;V)}^{1/2}\|\partial_t^2 u\|_{L^2(\R;V')}^{1/2}.
          \end{equation}
          Hence, we conclude that
          \[
            \widetilde{\widetilde{W}}(-\infty,\infty;V)=W^2(-\infty,\infty; V,V');
          \]
          see Definition~\ref{F:def: special function spaces}. Thus, we have proved the following result.
          
          \begin{corollary}
          \label{F:corollary: extension operator double tilde space}
              Let $T>0$ and suppose that $V\hookrightarrow H\hookrightarrow V'$ is a Gelfand triple. Then, for any $\Delta>0$, there exists a continuous extension operator 
            \[
                P^2\colon \widetilde{\widetilde{W}}(0,T;V)\to \widetilde{\widetilde{W}} (-\infty,\infty;V)
            \]
            such that $P^2(u)|_{(0,T)}=u$ and $\supp(P^2(u))\subset [-\Delta,T+\Delta]$, $\Delta>0$.
          \end{corollary}

    \subsection{Regularity of potential term for \texorpdfstring{$C^1(I;L^p(\Omega))$}{CLp} potentials}
    \label{F:appendix: proofs for energy identity}

    Here, we establish the proof of Claim~\ref{F:claim: derivative of potential term}. 
    
    \begin{proof}
        First, we show that for all Banach spaces $X\hookrightarrow Y$ with dense inclusion, any open, bounded interval $I\subset\R$, $1\leq p,q<\infty$ and any $w\in L^p(I;X)\cap W^{1,q}(I;Y)$ with mollification  
            \begin{equation}
            \label{F:eq: mollification}
                w_\delta\vcentcolon = \rho_\delta\ast W
            \end{equation}
            there holds:
            \begin{enumerate}[(a)]
                \item\label{F:prop 1 mollification} $w_\delta|_I\in C^{\infty}(I;X)$,
                \item\label{F:prop 2 mollification} $w_\delta\to w$ in $L^p(I;X)$,
                \item\label{F:prop 3 mollification}  $\partial_t w_\delta\to \partial_t w$ in $L^q(I;Y)$.
            \end{enumerate}
            Here, $W$ is the function obtained by first reflecting it at $t=a\vcentcolon =\inf I$ and $t=b\vcentcolon =\sup I$ via 
            \begin{equation}
           \label{F:eq: extension by reflection}
            \widetilde{w} (t)\vcentcolon =
            \begin{cases}
                w(t)&\text{ for }a\leq t\leq b,\\
                w(2a-t)&\text{ for }2a-b\leq t\leq a,\\
                w(2b-t)&\text{ for } b\leq t\leq 2b-a
            \end{cases}
        \end{equation}
            and then multiplying $\widetilde{w}$ by a cutoff function $\psi\in C_c^{\infty}$ such that $\psi=1$ in a neighborhood of $[0,T]$. It is easy to show that $\widetilde{w}\in L^p(2a-b,2b-a;X)$. Furthermore, by the density of $C^{\infty}(\overline{I};Y)$ in $W^{1,1}(I;Y)$ and the fundamental theorem of calculus, one can show that $\partial_t \widetilde{u}\in L^q(I;Y)$, which is given by
            \[
                \partial_t\widetilde{w} (t)\vcentcolon =
            \begin{cases}
                \partial_t w(t)&\text{ for }a\leq t\leq b,\\
                -\partial_t w(2a-t)&\text{ for }2a-b\leq t\leq a,\\
                -\partial_t w(2b-t)&\text{ for } b\leq t\leq 2b-a.
            \end{cases}
            \]
            Clearly, $W$ has the same properties. Exactly, as in the scalar valued case, one has
            \[
                \partial_t w_\delta=\rho_\delta\ast \partial_t W
            \]
            and, therefore, the assertions \ref{F:prop 1 mollification}--\ref{F:prop 3 mollification}   follow by standard results on the convolution.

            Hence, by the above preparation, we can choose sequences $u_\delta,v_\delta\in C^{\infty}(I;\widetilde{H}^s(\Omega))$ having the properties \ref{F:prop 1 mollification}--\ref{F:prop 3 mollification} with $X=\widetilde{H}^s(\Omega)$, $Y=\widetilde{L}^2(\Omega)$ and $p=q=2$. For such functions, we can calculate pointwise,
            \begin{equation}
            \label{F:eq: pointwise derivative}
                \partial_t \langle \mathcal{Q} u_\delta,v_\delta\rangle_{L^2(\Omega)}=\langle \partial_t u_\delta,\mathcal{Q} v_\delta\rangle_{L^2(\Omega)}+\langle \mathcal{Q} u_\delta,\partial_t v_\delta\rangle_{L^2(\Omega)}+\langle (\partial_t \mathcal{Q}) u_\delta,v_\delta\rangle_{L^2(\Omega)}.
            \end{equation}
            Here, we are using that $q(t)\colon \widetilde{H}^s(\Omega)\to \widetilde{L}^2(\Omega)$ is a continuous map for all $0\leq t\leq T$ (see Lemma \ref{F:lemma: cond on q}). Using the convergence properties of $u_\delta$ and $v_\delta$, $\mathcal{Q}\in C^1(\overline{I};L^p(\Omega))$ and Lemma \ref{F:lemma: cond on q}, we can pass to the limit $\delta\to 0$ and deduce formula \eqref{F:eq: derivative of potential term}. Finally, the regularity $W^{1,1}(I)$ follows from Lemma \ref{F:lemma: cond on q}.
    \end{proof}

    \subsection{Friedrich's lemma}
    \label{F:subsec: Friedrich's lemma}

    In this section we present Friedrich's lemma, which is discussed in \cite{lions2013equations} for the case of bounded linear operators $\mathcal{A}(t)$ from a Hilbert space $H$ to itself.

     \begin{lemma}[{Friedrich's lemma, cf. \cite[p. 72]{lions2013equations}}]
       \label{F:lemma: Friedrichs lemma}
           Let $V$ be a (complex) Hilbert space with (anti) duality pairing $\langle\cdot,\cdot\rangle$, $(\rho_\delta)_{\delta>0}$ the standard mollifier, and suppose that $\mathcal{A}(t)\in L(V,V')$, $t\in\R$, is a family of linear operators such that 
           \[
                t\mapsto \langle \mathcal{A}(t) u,v\rangle \in C^1(\R)
           \]
           for all $u,v\in V$ and 
            \begin{equation}
            \label{F:eq: bound on derivative Friedrich}
                |\partial_t \langle \mathcal{A}(t)u,v\rangle|=|\langle \partial_t \mathcal{A}(t) u,v\rangle|\leq C\|u\|_V \|v\|_V.
            \end{equation}
           Then, for all $u\in L^2(\R;V)$, there holds 
           \begin{equation}
           \label{F:eq: Friedrich convergence}
               \partial_t (\mathcal{A}(t)(\rho_\delta\ast u)-\rho_\delta\ast (\mathcal{A}(t)u))\to 0\text{ in }L^2(\R;V')
           \end{equation}
           as $\delta\to 0$.
       \end{lemma}

   \begin{proof}
Let $J\colon V\to V'$ denote the Riesz isomorphism,
characterized by
\[
    \langle Ju,v\rangle_{V'\times V}
    =\langle u,v\rangle_V,
    \qquad u,v\in V,
\]
and define
\[
    B(t)\coloneqq J^{-1}\mathcal A(t)\in\mathcal L(V).
\]
Then
\[
    \langle B(t)u,v\rangle_V
    =\langle\mathcal A(t)u,v\rangle_{V'\times V}
\]
for all $u,v\in V$. Consequently, the scalar functions
$t\mapsto\langle B(t)u,v\rangle_V$ belong to $C^1(\R)$ and satisfy
\[
    \left|\frac{d}{dt}\langle B(t)u,v\rangle_V\right|
    \leq C\|u\|_V\|v\|_V.
\]
Thus, $B$ satisfies the hypotheses of
\cite[Lemma~7.2, pp.~72--73]{lions2013equations}
with $\mathcal H=V$. For every $u\in L^2(\R;V)$,
that lemma yields
\[
    \partial_t\bigl[
        B(t)(\rho_\delta*u)-\rho_\delta*(Bu)
    \bigr]
    \to 0
    \quad\text{in }L^2(\R;V).
\]
Since $J$ is bounded and independent of time, it commutes
with convolution and distributional differentiation.
Applying $J$ to the preceding convergence therefore gives
\[
    \partial_t\bigl[
        \mathcal A(t)(\rho_\delta*u)
        -\rho_\delta*(\mathcal Au)
    \bigr]
    \to 0
    \quad\text{in }L^2(\R;V'),
\]
as claimed.
\end{proof}

    \section{Proof of the energy identity}
    \label{F:sec: proof of energy identity}

    Here we establish the energy identity Proposition~\ref{F:prop: energy identity weakly viscous nonlocal wave} by the method of double regularization.

     \begin{proof}[{Proof of Proposition~\ref{F:prop: energy identity weakly viscous nonlocal wave}}]
        We first prove the identity \eqref{F:eq: energy identity weakly viscous nonlocal wave eq} and then explain how to modify the argument to establish \eqref{F:eq: energy identity weakly viscous nonlocal wave eq continuous potential}. For any $0<t_0\leq T$, let us introduce the cutoff functions $\eta_\lambda\in C(\R)$, $0<\lambda<t_0$, by
        \begin{equation}
        \label{F:eq: cutoff}
            \eta_\lambda(t)=
            \begin{cases}
                 t/\lambda,&\text{ for }0\leq t\leq \lambda,\\
                1,&\text{ for }\lambda\leq t\leq t_0-\lambda,\\
                (t_0-t)/\lambda,&\text{ for }t_0-\lambda\leq t\leq t_0,\\
                0,&\text{ for } t\notin [0,t_0].
            \end{cases}
        \end{equation}
        For later convenience, we set $\eta_0\vcentcolon = \chi_{[0,t_0]}$. Moreover, let $(\rho_\delta)_{0<\delta<T/2}\subset C_c^{\infty}(\R)$ be the standard mollifiers in one space dimension.
        
        Let us recall that by Corollary~\ref{F:corollary: extension operator double tilde space} there exists a continuous extension operator 
        \begin{equation}
        \label{F:eq: cont extension m =2}
            P^2\colon \widetilde{\widetilde{W}}(0,T;\widetilde{H}^s(\Omega))\to \widetilde{\widetilde{W}}(-\infty,
            \infty;\widetilde{H}^s(\Omega));
        \end{equation}
        see Definition~\ref{F:def: special function spaces}.
        Furthermore, the construction of the extension operator $P^2$, carried out in Appendix~\ref{F:appendix: generlaized vector-valued Sobolev}, shows that if 
        \[
            u\in \widetilde{\widetilde{W}}(0,T;\widetilde{H}^s(\Omega))\text{ with }\partial_t u\in L^2(0,T;\widetilde{H}^r(\Omega)),
        \]
       then its extension $P^2(u)$ obeys
        \begin{equation}
        \label{F:eq: time regularity extension} 
            \partial_t P^2(u)\in L^2(\R;\widetilde{H}^r(\Omega)) \text{ with }\| \partial_t P^2(u)\|_{L^2(\R;\widetilde{H}^r(\Omega))}\leq C\|u\|_{H^1(0,T;\widetilde{H}^r(\Omega))},
        \end{equation}
        where the lower order term $\|u\|_{L^2(0,T;\widetilde H^r(\Omega))}$ comes from the cutoff used in the construction.
        Moreover, an inspection of the proof (see formulas \eqref{F:eq: def of coeff}, \eqref{F:eq: higher order reflection at t=0}, and \eqref{F:eq: formula for mth derivative}) demonstrates that if 
        \[
            u\in C_w([0,T];\widetilde{H}^s(\Omega))\text{ with }\partial_t u\in C_w([0,T];\widetilde{L}^2(\Omega))
        \]
        then 
        \[
            P^2(u)\in C_w(\R;\widetilde{H}^s(\Omega))\text{ with }\partial_t P^2(u)\in C_w(\R;\widetilde{L}^2(\Omega)).
        \]
        Similarly, if $q\in C^1([0,T];L^p(\Omega))$, then we can use the same extension procedure to get a compactly supported function $Q\in C^1(\R;L^p(\Omega))$ with $Q|_{[0,T]}=q$ and satisfying
        \begin{equation}
        \label{F:eq: estimate potential}
            \|Q\|_{C^1(\R;L^p(\Omega))}\leq C\|q\|_{C^1([0,T];L^p(\Omega))};
        \end{equation}
        see, for example, \eqref{F:eq: equal k der} and \eqref{F:eq: formula for mth derivative}.
        
        Thus, if we replace the weak finite energy solution $u$ by $U\vcentcolon = P^2(u)$ and $q$ by $Q$, then they have the following properties:
        \begin{enumerate}[(i)]
            \item\label{F:prop 1 for cont} $U\in L^2(\R;\widetilde{H}^s(\Omega))$ with $\partial_t U\in L^2(\R;\widetilde{H}^r(\Omega))$;
            \item\label{F:prop 2 for cont} $U\in C_w(\R;\widetilde{H}^s(\Omega))$ and $\partial_t U\in C_w(\R;\widetilde{L}^2(\Omega))$;
            \item\label{F:prop 3 for cont} $Q\in C^1(\R;L^p(\Omega))$;
            \item\label{F:support} $U$ and $Q$ have compact support;
            \item\label{F:prop 4 for cont}
            $U$ satisfies \eqref{F:eq: existence weakly viscous case} in the sense of distributions on $\Omega_T$ and $(U,\partial_t U)(0)=(u_0,u_1)$ in $\widetilde{H}^s(\Omega)\times L^2(\Omega)$ (see Lemma~\ref{F:lemma: u1 is attained strongly} and \eqref{F:eq: equal k der}).
        \end{enumerate}
        Furthermore, we extend $f$ outside of the interval $(0,T)$ by zero to get a function $F\in L^2(\R;H^{-r}(\Omega))$. For the subsequent arguments, we need the following result, whose proof is established in Section~\ref{F:appendix: proofs for energy identity}.
        \begin{claim}
        \label{F:claim: derivative of potential term}
            Let $I\subset\R$ be any open and bounded interval, $s>0$, and assume that $2\leq p\leq\infty$ satisfies \eqref{F:eq: exponent p}. Then, for any  real-valued potential $\mathcal{Q}\in C^1(\overline I;L^p(\Omega))$ and all functions $u,v\in L^2(I;\widetilde{H}^s(\Omega))$ with $\partial_t u,\partial_t v\in L^2(I;\widetilde{L}^2(\Omega))$, we have 
            \begin{equation}
            \label{F:eq: regularity of potential term}
                I\ni t\mapsto \langle \mathcal{Q}(t) u(t),v(t)\rangle_{L^2(\Omega)} \in W^{1,1}(I)
            \end{equation}
            and there holds
            \begin{equation}
            \label{F:eq: derivative of potential term}
                \partial_t \langle \mathcal{Q} u,v\rangle_{L^2(\Omega)}=\langle \partial_t u,\mathcal{Q} v\rangle_{L^2(\Omega)}+\langle \mathcal{Q} u,\partial_t v\rangle_{L^2(\Omega)}+\langle (\partial_t \mathcal{Q}) u,v\rangle_{L^2(\Omega)}.
            \end{equation}
        \end{claim}
        The rest of the proof is divided into four steps. \\
        
       \noindent \textit{Step 1.} \textit{Control of the elliptic operator $A(t)=(-\Delta)^s+Q(t)$.} Using Claim  \ref{F:claim: derivative of potential term}, the properties \ref{F:prop 1 for cont}--\ref{F:prop 4 for cont}, and the fundamental theorem of calculus, we deduce that
       \[
       \begin{split}
            0&=\int_{\R}\partial_t \langle A(\rho_\delta\ast \eta_\lambda U),(\rho_\delta\ast \eta_\lambda U)\rangle\,dt\\
            &=\int_{\R}\langle (\partial_tA) (\rho_\delta\ast \eta_\lambda U),(\rho_\delta\ast \eta_\lambda U)\rangle\,dt\\
            &\quad +2\re\int_{\R}\langle A (\rho_\delta\ast \eta_\lambda U),\partial_t (\rho_\delta\ast \eta_\lambda U)\rangle\,dt.
       \end{split}
       \]
       Here and in the rest of this proof, we simply write $\langle\cdot,\cdot\rangle$ to denote the pairing between $\widetilde{H}^s(\Omega)$ and $H^{-s}(\Omega)$ or the inner product in $\widetilde{L}^2(\Omega)$. Using
       \[
       \begin{split}
            \partial_t (\rho_\delta\ast \eta_\lambda U)&=\rho'_\delta\ast \eta_\lambda U,\\
             \partial_t (\rho_\delta\ast \eta_\lambda U)&=\rho_\delta\ast \eta'_\lambda U+\rho_\delta\ast \eta_\lambda (\partial_t U),
       \end{split}
       \]
       we obtain
       \begin{equation}
       \label{F:eq: identity for limi in lambda}
            \begin{split}
                0&= \int_{\R}\langle (\partial_tQ) (\rho_\delta\ast \eta_\lambda U),(\rho_\delta\ast \eta_\lambda U)\rangle\,dt\\
            &\quad +2\re \int_{\R}\langle  \rho_\delta\ast(A\eta_\lambda U),\partial_t (\rho_\delta\ast \eta_\lambda U)\rangle\,dt\\
            &\quad +2\re \int_{\R}\langle  A(\rho_\delta\ast\eta_\lambda U)-\rho_\delta\ast (A\eta_\lambda U),\partial_t (\rho_\delta\ast \eta_\lambda U)\rangle\,dt\\
            &=\int_{\R}\langle (\partial_tQ) (\rho_\delta\ast \eta_\lambda U),(\rho_\delta\ast \eta_\lambda U)\rangle\,dt\\
             &\quad +2\re \int_{\R}\langle  \rho_\delta\ast(A\eta_\lambda U),\rho_\delta\ast \eta'_\lambda U\rangle\,dt+2\re \int_{\R}\langle  \rho_\delta\ast(A\eta_\lambda U),\rho_\delta\ast \eta_\lambda (\partial_t U)\rangle\,dt\\
             &\quad +2\re \int_{\R}\langle  A(\rho_\delta\ast\eta_\lambda U)-\rho_\delta\ast (A\eta_\lambda U),\rho'_\delta\ast \eta_\lambda U\rangle\,dt
            \end{split}
       \end{equation}
       Now, we wish to pass to the limit $\lambda\to 0$. We clearly have
       \begin{enumerate}[{(N}1)]
           \item\label{F:prop 1 eta lambda} $|\eta_\lambda(t)|\leq 1$,
           \item\label{F:prop 2 eta lambda} $\eta'_\lambda(t)=\lambda^{-1}(\chi_{(0,\lambda)}-\chi_{(t_0-\lambda,t_0)})$ with $\int_{\R}|\eta'_\lambda|\,dt=2$
           \item\label{F:prop 3 eta lambda} and $\eta_\lambda(t)\to \eta_0(t)\text{ as }\lambda\to 0$ for a.e. $t$.
       \end{enumerate}
       Using Lebesgue's dominated convergence theorem, the properties of $U$ and $\mathcal{Q}\in C^1_c(\R;L^p(\Omega))$, we infer that
       \begin{equation}
       \label{F:eq: conv eta lambda times U}
       \begin{cases}
           \eta_\lambda U\to \eta_0 U\text{ in }L^2(\R;\widetilde{H}^s(\Omega)),\\
           \eta_\lambda \partial_t U\to \eta_0 \partial_t U\text{ in }L^2(\R;\widetilde{H}^r(\Omega)),\\
            A\eta_\lambda U\to  A\eta_0 U \text{ in }L^{2}(\R;H^{-s}(\Omega))
       \end{cases}
       \end{equation}
       as $\lambda\to 0$. From Young's inequality, \eqref{F:eq: conv eta lambda times U} and $\mathcal{Q}\in C^1_c(\R;L^p(\Omega))$, we deduce
       \begin{equation}
       \label{F:eq: convergence 1 of int}
           \begin{split}
               &\int_{\R}\langle (\partial_tQ) (\rho_\delta\ast \eta_\lambda U),(\rho_\delta\ast \eta_\lambda U)\rangle\,dt\to \int_{\R}\langle (\partial_tQ) (\rho_\delta\ast \eta_0 U),(\rho_\delta\ast \eta_0 U)\rangle\,dt,\\
               &\re \int_{\R}\langle  A(\rho_\delta\ast\eta_\lambda U)-\rho_\delta\ast (A\eta_\lambda U),\rho'_\delta\ast \eta_\lambda U\rangle\,dt\\
               &\to \re \int_{\R}\langle  A(\rho_\delta\ast\eta_0 U)-\rho_\delta\ast (A\eta_0 U),\rho'_\delta\ast \eta_0 U\rangle\,dt
           \end{split}
       \end{equation}
       as $\lambda\to 0$. Next, we assert that there holds
       \begin{equation}
       \label{F:eq: convergence 2 of int}
           \re \int_{\R}\langle  \rho_\delta\ast(A\eta_\lambda U),\rho_\delta\ast \eta_\lambda (\partial_t U)\rangle\,dt\to \re \int_{\R}\langle  \rho_\delta\ast(A\eta_0 U),\rho_\delta\ast \eta_0(\partial_t U)\rangle\,dt.
       \end{equation}
       First, let us note that we have
       \begin{equation}
       \label{F:eq: regularity of mollified time-derivative}
        \rho_\delta\ast \eta_\lambda \partial_t U (t)=\rho'_\delta\ast \eta_\lambda U(t)-\rho_\delta\ast \eta'_\lambda U(t)\in \widetilde{H}^s(\Omega)
       \end{equation}
       for $t\in\R$. Then, using Young's inequality, \ref{F:prop 1 eta lambda}, and \ref{F:prop 2 eta lambda}, we can estimate
       \[
       \begin{split}
           &\|\rho_\delta\ast \eta_\lambda \partial_t U (t)\|_{\widetilde{H}^s(\Omega)}\leq \|\rho'_\delta\|_{L^2(\R)}\|\eta_\lambda U\|_{L^2(\R;\widetilde{H}^s(\Omega))}+\|\rho_\delta\|_{L^{\infty}(\R)}\|\eta'_\lambda U\|_{L^1(\R;\widetilde{H}^s(\Omega))}\\
           &\leq \|\rho'_\delta\|_{L^2(\R)}\|\eta_\lambda\|_{L^{\infty}(\R)}\| U\|_{L^2(\R;\widetilde{H}^s(\Omega))}+\|\rho_\delta\|_{L^{\infty}(\R)}\|\eta'_\lambda\|_{L^1(\R)} \|U\|_{L^{\infty}(\R;\widetilde{H}^s(\Omega))}\\
           &\leq C\|\rho'_\delta\|_{L^2(\R)}\|U\|_{L^2(\R;\widetilde{H}^s(\Omega))}+\|\rho_\delta\|_{L^{\infty}(\R)}\|U\|_{L^{\infty}(\R;\widetilde{H}^s(\Omega))}.
       \end{split}
       \]
       Hence, $(\rho_\delta\ast\eta_\lambda \partial_t U)_{0<\lambda<t_0}$ is uniformly bounded in $L^{\infty}(\R;\widetilde{H}^s(\Omega))$ and thus there exists $\chi_\delta\in L^{\infty}(\R;\widetilde{H}^s(\Omega))$ such that
       \begin{equation}
       \label{F:eq: mollification time der eta null}
            \rho_\delta\ast\eta_\lambda \partial_t U\weakstar \chi_\delta\text{ in }L^{\infty}(\R;\widetilde{H}^s(\Omega))
       \end{equation}
       as $\lambda\to 0$. On the other hand, by \eqref{F:eq: conv eta lambda times U} and Young's inequality, we have
       \[
        \rho_\delta\ast \eta_\lambda \partial_t U\to \rho_\delta\ast \eta_0 \partial_t U\text{ in }L^{\infty}(\R;\widetilde{L}^2(\Omega)).
       \]
       Therefore, we may conclude that
       \begin{equation}
        \label{F:eq: mollification time der eta null 2}
           \chi_\delta=\rho_\delta\ast \eta_0 \partial_t U
       \end{equation}
        and 
       \begin{equation}
       \label{F:eq: weak star convergence of mollified time der}
            \rho_\delta\ast\eta_\lambda \partial_t U\weakstar \chi_\delta\text{ in }L^{\infty}(\R;\widetilde{H}^s(\Omega))
       \end{equation}
       as $\lambda\to 0$. Furthermore, we can write
        \begin{equation}
        \label{F:eq: decomposition of AU}
            \rho_\delta\ast \eta_\lambda AU=\rho_\delta\ast (\eta_\lambda-\eta_0)AU+\rho_\delta\ast \eta_0 AU.
        \end{equation}
        Using \eqref{F:eq: conv eta lambda times U}, the compact support assumption of $U$ and Young's inequality, we see that the first term in  \eqref{F:eq: decomposition of AU} goes to zero in $L^1(\R;H^{-s}(\Omega))$ as $\lambda\to 0$ and the second term in \eqref{F:eq: decomposition of AU} belongs to $L^1(\R;H^{-s}(\Omega))$. Thus, using these properties and \eqref{F:eq: weak star convergence of mollified time der}, we deduce the convergence \eqref{F:eq: convergence 2 of int}. 

       The remaining term in \eqref{F:eq: identity for limi in lambda} can be written as
       \begin{equation}
       \label{F:eq: reform remaining term}
       \begin{split}
           &\re \int_{\R}\langle  \rho_\delta\ast(A\eta_\lambda U),\rho_\delta\ast \eta'_\lambda U\rangle\,dt\\
           &=\re \int_{\R}\langle  \rho_\delta\ast(A(\eta_\lambda-\eta_0) U),\rho_\delta\ast \eta'_\lambda U\rangle\,dt\\
           &\quad +\re \int_{\R}\langle  \rho_\delta\ast(A\eta_0 U),\rho_\delta\ast \eta'_\lambda U\rangle\,dt.
        \end{split}
       \end{equation}
       Moreover, \ref{F:prop 2 eta lambda} implies  that the first term on the right-hand side of \eqref{F:eq: reform remaining term} goes to zero as $\lambda\to 0$. More precisely, using Young's inequality we can estimate
       \[
       \begin{split}
            &\left|\re \int_{\R}\langle  \rho_\delta\ast(A(\eta_\lambda-\eta_0) U),\rho_\delta\ast \eta'_\lambda U\rangle\,dt\right|\\
            &\quad \leq \|\rho_\delta\ast(A(\eta_\lambda-\eta_0) U) \|_{L^{\infty}(\R;H^{-s}(\Omega))}\|\rho_\delta\ast \eta'_\lambda U\|_{L^1(\R;\widetilde{H}^s(\Omega))}\\
            &\lesssim \|A(\eta_\lambda-\eta_0) U \|_{L^2(\R;H^{-s}(\Omega))}\| \eta'_\lambda U\|_{L^1(\R;\widetilde{H}^s(\Omega))}\\
            &\lesssim\|\eta'_\lambda\|_{L^1(\R)}\|(\eta_\lambda-\eta_0) U \|_{L^2(\R;\widetilde{H}^s(\Omega))}\| U\|_{L^{\infty}(\R;\widetilde{H}^s(\Omega))}\\
            &\lesssim\|(\eta_\lambda-\eta_0) U \|_{L^2(\R;\widetilde{H}^s(\Omega))}\to 0
       \end{split}
       \]
       as $\lambda\to 0$. In the last inequality, we used that $U\in C_w(\R; \widetilde{H}^s(\Omega))$, \ref{F:prop 2 eta lambda} and \eqref{F:eq: conv eta lambda times U}. The second term in \eqref{F:eq: reform remaining term}, can be written as
       \[
       \begin{split}
           &\int_{\R}\langle  \rho_\delta\ast(A\eta_0 U),\rho_\delta\ast \eta'_\lambda U\rangle\,dt=\int_{\R}\langle  \rho_\delta\ast \rho_\delta\ast(A\eta_0 U),\eta'_\lambda U\rangle\,dt\\
           &\quad =\fint_0^\lambda \langle  \rho_\delta\ast \rho_\delta\ast(A\eta_0 U),\eta_0 U\rangle\,dt-\fint_{t_0-\lambda}^{t_0} \langle  \rho_\delta\ast \rho_\delta\ast(A\eta_0 U),\eta_0 U\rangle\,dt.
       \end{split}
       \]
       Here, we use that $\rho_\delta$ is an even function. Next, note that the function
       \[
        t\mapsto  \rho_\delta\ast \rho_\delta\ast(A\eta_0 U)\in H^{-s}(\Omega)
       \]
       is continuous as $A\eta_0 U\in L^2(\R;H^{-s}(\Omega))$. Taking into account that the solution $U$ belongs to the space $C_w(\R;\widetilde{H}^s(\Omega))$, we conclude that 
       \[
       t\mapsto \langle  \rho_\delta\ast \rho_\delta\ast(A\eta_0 U), U\rangle
       \]
       is continuous and hence Lebesgue's differentiation theorem ensures that 
       \begin{equation}
       \label{F:eq: convergence last term}
       \begin{split}
            \int_{\R}\langle  \rho_\delta\ast(A\eta_0 U),\rho_\delta\ast \eta'_\lambda U\rangle\,dt&\to \sum_{j=0,1}(-1)^j\langle  \rho_\delta\ast \rho_\delta\ast(A\eta_0 U),\eta_0 U\rangle (j t_0)
       \end{split}
       \end{equation}
       as $\lambda\to 0$. Thus, the convergence results \eqref{F:eq: identity for limi in lambda}, \eqref{F:eq: convergence 1 of int}, \eqref{F:eq: convergence 2 of int}, and \eqref{F:eq: convergence last term} yield
       \begin{equation}
       \label{F:eq: identity for elliptic part}
           \begin{split}
               0&=\int_{\R}\langle (\partial_tQ) (\rho_\delta\ast \eta_0 U),(\rho_\delta\ast \eta_0 U)\rangle\,dt\\
               &+2\re \int_{\R}\langle  \rho_\delta\ast(A\eta_0 U),\rho_\delta\ast \eta_0(\partial_t U)\rangle\,dt\\
               &+2\re \int_{\R}\langle  A(\rho_\delta\ast\eta_0 U)-\rho_\delta\ast (A\eta_0 U),\rho'_\delta\ast \eta_0 U\rangle\,dt\\
               &+2\re\sum_{j=0,1}(-1)^j\langle  \rho_\delta\ast \rho_\delta\ast(A\eta_0 U),\eta_0 U\rangle (j t_0).
           \end{split}
       \end{equation}
       \noindent \textit{Step 2.} \textit{Control of the second order derivative.} From the support assumption on $U$ we infer
       \[
       \begin{split}
           0&=\int_\R \partial_t \langle \rho_\delta\ast \eta_\lambda \partial_t U, \rho_\delta\ast \eta_\lambda \partial_t U\rangle\,dt\\
           &=2\re\int_\R \langle \rho_\delta\ast \eta_\lambda \partial_t U, \rho_\delta\ast \eta_\lambda \partial_t^2 U\rangle\,dt\\
           &+2\re\int_\R \langle \rho_\delta\ast \eta_\lambda \partial_t U, \rho_\delta\ast \eta'_\lambda \partial_t U\rangle\,dt\\
           &=2\re\int_\R \langle \rho_\delta\ast \eta_\lambda \partial_t U, \rho_\delta\ast \eta_\lambda \partial_t^2 U\rangle\,dt\\
           &+2\re\int_\R \langle \rho_\delta\ast (\eta_\lambda-\eta_0) \partial_t U, \rho_\delta\ast \eta'_\lambda \partial_t U\rangle\,dt\\
           &+2\re\int_\R \langle \rho_\delta\ast \eta_0 \partial_t U, \rho_\delta\ast \eta'_\lambda \partial_t U\rangle\,dt.
       \end{split}
       \]
       Arguing as in Step 1, we see that
       \begin{equation}
           \begin{split}
               &\re\int_\R \langle \rho_\delta\ast (\eta_\lambda-\eta_0) \partial_t U, \rho_\delta\ast \eta'_\lambda \partial_t U\rangle\,dt\to 0,\\
               &\re\int_\R \langle \rho_\delta\ast \eta_0 \partial_t U, \rho_\delta\ast \eta'_\lambda \partial_t U\rangle\,dt\to\re \sum_{j=0,1}(-1)^j\langle\rho_\delta\ast  \rho_\delta\ast \eta_0 \partial_t U, \eta_0 \partial_t U\rangle(j t_0)
           \end{split}
       \end{equation}
       as $\lambda\to 0$. In fact, for the first term, we use \ref{F:prop 2 eta lambda}, \eqref{F:eq: conv eta lambda times U}, $\partial_t U\in L^{\infty}(\R;\widetilde{L}^2(\Omega))$, Young's inequality and H\"older's inequality to estimate
       \begin{equation}
       \label{F:eq: estimate for eta prime}
       \begin{split}
           &\left|\re\int_\R \langle \rho_\delta\ast (\eta_\lambda-\eta_0) \partial_t U, \rho_\delta\ast \eta'_\lambda \partial_t U\rangle\,dt\right|\\
           &\leq \|\rho_\delta\ast (\eta_\lambda-\eta_0) \partial_t U\|_{L^{\infty}(\R;L^2(\Omega))}\|\rho_\delta\ast \eta'_\lambda \partial_t U\|_{L^1(\R;L^2(\Omega))}\\
           &\lesssim\|(\eta_\lambda-\eta_0) \partial_t U\|_{L^{2}(\R;L^2(\Omega))}\|\eta'_\lambda \partial_t U\|_{L^1(\R;L^2(\Omega))}\\
           &\lesssim\|(\eta_\lambda-\eta_0) \partial_t U\|_{L^{2}(\R;L^2(\Omega))}\|\eta'_\lambda\|_{L^1(\R)}\| \partial_t U\|_{L^{\infty}(\R;L^2(\Omega))}\\
           &\lesssim\|(\eta_\lambda-\eta_0) \partial_t U\|_{L^{2}(\R;L^2(\Omega))}\| \partial_t U\|_{L^{\infty}(\R;L^2(\Omega))}\to 0
       \end{split}
       \end{equation}
       as $\lambda\to 0$. For the second term, we note that \ref{F:prop 2 eta lambda} implies
       \[
       \begin{split}
            &\re\int_\R \langle \rho_\delta\ast \eta_0 \partial_t U, \rho_\delta\ast \eta'_\lambda \partial_t U\rangle\,dt=\re\int_\R \langle \rho_\delta\ast\rho_\delta\ast \eta_0 \partial_t U, \eta'_\lambda \partial_t U\rangle\,dt\\
            &=\re\fint_0^{\lambda} \langle \rho_\delta\ast\rho_\delta\ast \eta_0 \partial_t U, \eta_0\partial_t U\rangle\,dt-\re\fint_{t_0-\lambda}^{t_0} \langle \rho_\delta\ast\rho_\delta\ast \eta_0 \partial_t U, \eta_0\partial_t U\rangle\,dt
       \end{split}
       \]
        and we use that the mapping
       \[
            t\mapsto \langle \rho_\delta\ast \rho_\delta\ast \eta_0 \partial_t U, \eta_0 \partial_t U\rangle
       \]
       is continuous, since $\partial_t U\in C_w(\R;\widetilde{L}^2(\Omega))$ and the first factor is a continuous function to $\widetilde{L}^2(\Omega)$, to conclude via Lebesgue's differentiation theorem that
       it converges to 
       \[
        \re\sum_{j=0,1}(-1)^j\langle\rho_\delta\ast  \rho_\delta\ast \eta_0 \partial_t U,  \eta_0\partial_t U\rangle(j t_0).
       \]
       On the other hand, in the limit $\lambda\to 0$ we have
       \begin{equation}
       \label{F:eq: convergence 3 of int}
        \re\int_\R \langle \rho_\delta\ast \eta_\lambda \partial_t U, \rho_\delta\ast \eta_\lambda \partial_t^2 U\rangle\,dt\to \re\int_\R \langle \rho_\delta\ast \eta_0 \partial_t U, \rho_\delta\ast \eta_0 \partial_t^2 U\rangle\,dt.
       \end{equation}
        The proof of the identity \eqref{F:eq: convergence 3 of int} is completely analogous to \eqref{F:eq: convergence 2 of int}.
       
       Therefore, we have established that there holds
       \begin{equation}
       \label{F:eq: identity for second order derivative}
           \begin{split}
                &\re\int_\R \langle \rho_\delta\ast \eta_0 \partial_t U, \rho_\delta\ast \eta_0 \partial_t^2 U\rangle\,dt\\
                &+\re\sum_{j=0,1}(-1)^j\langle \rho_\delta\ast \rho_\delta\ast \eta_0 \partial_t U,  \eta_0\partial_t U\rangle(j t_0)=0.
           \end{split}
       \end{equation}
       \noindent \textit{Step 3.} \textit{Approximate energy identity.} By adding \eqref{F:eq: identity for elliptic part} to twice \eqref{F:eq: identity for second order derivative}, we get
       \[
       \begin{split}
           0&=\int_{\R}\langle (\partial_tQ) (\rho_\delta\ast \eta_0 U),(\rho_\delta\ast \eta_0 U)\rangle\,dt\\
               &+2\re \int_{\R}\langle  \rho_\delta\ast \eta_0(\partial_t^2+A) U,\rho_\delta\ast \eta_0(\partial_t U)\rangle\,dt\\
               &+2\re \int_{\R}\langle  A(\rho_\delta\ast\eta_0 U)-\rho_\delta\ast (A\eta_0 U),\rho'_\delta\ast \eta_0 U\rangle\,dt\\
               &+2\re\sum_{j=0,1}(-1)^j\langle \rho_\delta\ast \rho_\delta\ast \eta_0 \partial_t U, \eta_0\partial_t U\rangle(j t_0)\\
               &+2\re\sum_{j=0,1}(-1)^j\langle  \rho_\delta\ast \rho_\delta\ast(A\eta_0 U),\eta_0 U\rangle (j t_0).
       \end{split}
       \]
       Next, we take into account that $U$ solves \eqref{F:eq: existence weakly viscous case} in $\Omega_T$ and therefore by approximation $U$ satisfies
       \begin{equation}
       \label{F:eq: weak form PDE}
        \left(\partial_t^2+(-\Delta)^r \partial_t+A\right)U=F\quad\text{in}\quad L^2(0,T;H^{-s}(\Omega)).
       \end{equation}
       Hence, using $\rho_\delta\ast \eta_0 U\in L^2(\R;\widetilde{H}^s(\Omega))$ and $\rho_\delta\ast \eta_0 \partial_t U\in L^{\infty}(\R;\widetilde{H}^s(\Omega))$ (see \eqref{F:eq: mollification time der eta null}--\eqref{F:eq: mollification time der eta null 2}), we get
       \begin{equation}
       \label{F:eq: identity to pass to the lmit delta to 0}
        \begin{split}
           0&=\int_{\R}\langle (\partial_tQ) (\rho_\delta\ast \eta_0 U),(\rho_\delta\ast \eta_0 U)\rangle\,dt\\
               &-2\re \int_{\R}\langle  \rho_\delta\ast \eta_0(-\Delta)^{r/2} \partial_t U,\rho_\delta\ast \eta_0 (-\Delta)^{r/2} \partial_t U\rangle\,dt\\
               &+2\re \int_{\R}\langle  \rho_\delta\ast F\eta_0 ,\rho_\delta\ast \eta_0(\partial_t U)\rangle\,dt\\
               &-2\re \int_{\R}\langle  \partial_t(A(\rho_\delta\ast\eta_0 U)-\rho_\delta\ast (A\eta_0 U)),\rho_\delta\ast \eta_0 U\rangle\,dt\\
               &+2\re\sum_{j=0,1}(-1)^j\langle \rho_\delta\ast \rho_\delta\ast \eta_0 \partial_t U, \eta_0\partial_t U\rangle(j t_0)\\
               &+2\re\sum_{j=0,1}(-1)^j\langle  \rho_\delta\ast \rho_\delta\ast(A\eta_0 U),\eta_0 U\rangle (j t_0).
       \end{split}
       \end{equation}
       For the term involving the commutator, we first move $\rho'_\delta$ to the first factor, which gives a minus sign as $\rho'_\delta$ is odd, then let $\partial_t$ act on the commutator, and afterwards move back $\rho_\delta$ to the second factor in the antiduality pairing.\\
       
       \noindent \textit{Step 4.} \textit{Limit $\delta\to 0$.} Finally, we wish to pass to the limit $\delta\to 0$ in \eqref{F:eq: identity to pass to the lmit delta to 0}. Since $\partial_t Q\in C_c(\R;L^p(\Omega))$ induces a continuous map from $L^2(\R;\widetilde{H}^s(\Omega))$ to $L^2(\R;\widetilde{L}^2(\Omega))$ (see Lemma~\ref{F:lemma: cond on q}), we deduce from the properties of mollification that
       \begin{equation}
       \label{F:eq: potential term delta to zero}
           \int_{\R}\langle (\partial_tQ) (\rho_\delta\ast \eta_0 U),(\rho_\delta\ast \eta_0 U)\rangle\,dt\to \int_{0}^{t_0}\langle (\partial_tQ) U, U\rangle\,dt
       \end{equation}
       as $\delta\to 0$. Similarly, by properties of mollification and $\partial_t U\in L^2(\R;\widetilde{H}^r(\Omega))$, we obtain in the limit $\delta\to 0$:
       \begin{equation}
       \label{F:eq: viscous term delta to zero}
           \re \int_{\R}\langle  \rho_\delta\ast \eta_0(-\Delta)^{r/2} \partial_t U,\rho_\delta\ast \eta_0 (-\Delta)^{r/2} \partial_t U\rangle\,dt\to \int_{0}^{t_0}\| (-\Delta)^{r/2} \partial_t U\|_{L^2(\R^n)}^2\,dt.
       \end{equation}
       Using that $F\in L^2(\R;H^{-r}(\Omega))$ and the same reasoning as above, we have
        \begin{equation}
       \label{F:eq: source term delta to zero}
           \re \int_{\R}\langle  \rho_\delta\ast F\eta_0 ,\rho_\delta\ast \eta_0(\partial_t U)\rangle\,dt\to \re \int_{0}^{t_0}\langle   F ,\partial_t U\rangle\,dt
       \end{equation}
       for $\delta\to 0$.
       For the fourth term in \eqref{F:eq: identity to pass to the lmit delta to 0} goes to zero by Friedrich's lemma; see Appendix~\ref{F:subsec: Friedrich's lemma}. 
       
        For the boundary terms, let us introduce
       \begin{equation}
       \label{F:eq: gamma}
           \gamma\vcentcolon = \rho\ast \rho\text{ and } \gamma_\delta \vcentcolon = \rho_\delta\ast\rho_\delta
       \end{equation}
       for $0<\delta<1$. Observe that $\gamma$ and 
       $\gamma_\delta$, $0<\delta<1$, have the following properties:
       \begin{enumerate}[(I)]
           \item\label{F:prop 1 sigma} $\gamma,\gamma_\delta\in C_c^{\infty}(\R)$ are nonnegative, even functions,
           \item\label{F:prop 2 sigma} $\gamma_\delta=\delta^{-1}\gamma(t/\delta)$,
           \item\label{F:prop 3 sigma} $\gamma_\delta\to \delta_0$ in the sense of distributions,
           \item\label{F:prop 4 sigma} and, for all $0<\delta<t_0/2$, we have
           \begin{equation}
           \label{F:eq: integral sigma delta}
                \int_0^{t_0}\gamma_\delta\,dt=\int_{-t_0}^0 \gamma_\delta\,dt=1/2.
           \end{equation}
       \end{enumerate}
       So, using \eqref{F:eq: integral sigma delta}, we deduce that
       \begin{equation}
       \label{F:eq: bdry term delta to zero}
       \begin{split}
           &\re \langle \gamma_\delta\ast \eta_0 \partial_t U, \eta_0\partial_t U\rangle (t_1)= \re\int_{\R}\gamma_\delta (t) \langle (\eta_0\partial_t U)(t_1-t),(\eta_0\partial_t U)(t_1)\rangle \,dt\\
           &=\re\int_{t_1-t_0}^{t_1}\gamma_\delta (t)  \langle \partial_t U(t_1-t),\partial_t U(t_1)\rangle \,dt\\
           &=\re\int_{t_1-t_0}^{t_1}\gamma_\delta (t) \langle(\partial_t U(t_1-t)-\partial_t U(t_1)),\partial_t U(t_1)\rangle \,dt\\
           &\quad+\re\int_{t_1-t_0}^{t_1}\gamma_\delta (t) \langle\partial_t U(t_1),\partial_t U(t_1)\rangle \,dt\\
           &=\re\int_{t_1-t_0}^{t_1}\gamma_\delta (t) \langle(\partial_t U(t_1-t)-\partial_t U(t_1)),\partial_t U(t_1)\rangle \,dt+\frac{1}{2}\|\partial_t U(t_1)\|^2_{L^2(\Omega)}
       \end{split}
       \end{equation}
       for $t_1\in \{0,t_0\}$ and a.e. $t_0$. For the above computation, recall that by Lemma~\ref{F:lemma: u1 is attained strongly} we have $\partial_t U(0)=u_1\in \widetilde{L}^2(\Omega)$. Since $\partial_t U\in C_w(\R;\widetilde{L}^2(\Omega))$, we know that 
       \[
        t\mapsto \langle(\partial_t U(t_1-t)-\partial_t U(t_1)),\partial_t U(t_1)\rangle 
       \]
       is continuous and thus by \ref{F:prop 3 sigma} the second last term in \eqref{F:eq: bdry term delta to zero} converges to zero as $\delta\to 0$,
       \begin{equation}
       \label{F:eq: remainder 1 intial cond goes to zero}
           \re\int_{t_1-t_0}^{t_1}\gamma_\delta (t) \langle(\partial_t U(t_1-t)-\partial_t U(t_1)),\partial_t U(t_1)\rangle \,dt\to 0.
       \end{equation}
       Similarly, we can show the convergence of the last term in \eqref{F:eq: identity to pass to the lmit delta to 0}. In fact, as above, we have
       \begin{equation}
       \label{F:eq: convergence second initial cond}
       \begin{split}
           &\re\langle  \gamma_\delta\ast(A\eta_0 U),\eta_0 U\rangle (t_1)=\re\int_{\R}\gamma_\delta (t) \langle (A\eta_0 U)(t_1-t),(\eta_0 U)(t_1)\rangle \,dt\\
           &=\re\int_{t_1-t_0}^{t_1}\gamma_\delta (t)  \langle (A U)(t_1-t),U(t_1)\rangle \,dt\\
           &=\re\int_{t_1-t_0}^{t_1}\gamma_\delta (t) \langle((A U)(t_1-t)-(A U)(t_1)),U(t_1)\rangle \,dt\\
           &\quad +\re\int_{t_1-t_0}^{t_1}\gamma_\delta (t) \langle(A U)(t_1), U(t_1)\rangle \,dt\\
           &=\re\int_{t_1-t_0}^{t_1}\gamma_\delta (t) \langle((A U)(t_1-t)-(A U)(t_1)), U(t_1)\rangle \,dt\\
           &\quad+\frac{1}{2}\left(\|(-\Delta)^{s/2}U(t_1)\|^2_{L^2(\R^n)}+\langle Q(t_1)U(t_1),U(t_1)\rangle_{L^2(\Omega)}\right)
       \end{split}
       \end{equation}
       for $t_1\in \{0,t_0\}$ and a.e. $t_0$. Here, we also apply that the solution $U$ satisfies $U(0)=u_0\in \widetilde{H}^s(\Omega)$ (see Lemma~\ref{F:lemma: u1 is attained strongly}). Using $U\in C_w(\R;\widetilde{H}^s (\Omega))$ and Claim  \ref{F:claim: derivative of potential term}, we see that
       \[
        t\mapsto \langle((A U)(t_1-t)-(A U)(t_1)), U(t_1)\rangle
       \]
       is continuous. Hence, \ref{F:prop 3 sigma} ensures that
       \begin{equation}
       \label{F:eq: remainder 2 intial cond goes to zero}
        \re\int_{t_1-t_0}^{t_1}\gamma_\delta (t) \langle((A U)(t_1-t)-(A U)(t_1)), U(t_1)\rangle \,dt\to 0
       \end{equation}
       as $\delta\to 0$. Combining the above results, we get the energy identity:
       \begin{equation}
        \begin{split}
            0&= \int_{0}^{t_0}\langle (\partial_tQ) U, U\rangle\,dt-2\int_{0}^{t_0}\| (-\Delta)^{r/2} \partial_t U\|_{L^2(\R^n)}^2\,dt+2\re \int_{0}^{t_0}\langle   F ,\partial_t U\rangle\,dt\\
             &+\sum_{j=0,1}(-1)^j\left(\|\partial_t U(jt_0)\|^2_{L^2(\Omega)}+\|(-\Delta)^{s/2}U(jt_0)\|^2_{L^2(\R^n)}\right)\\
             &+\sum_{j=0,1}(-1)^j\langle Q(jt_0)U(jt_0),U(jt_0)\rangle_{L^2(\Omega)}.
        \end{split}
       \end{equation}   
       Using that $(U,\partial_t U)(0)=(u_0,u_1)$ in $\widetilde{H}^s(\Omega)\times L^2(\Omega)$ we can conclude the proof of the energy identity \eqref{F:eq: energy identity weakly viscous nonlocal wave eq}. \\

       Finally, we show how the above proof can be modified to obtain \eqref{F:eq: energy identity weakly viscous nonlocal wave eq continuous potential} for potentials $q$ only satisfying the conditions of Lemma \ref{F:lemma: cond on q}. As above, we extend $u$, $f$ and $q$ to functions $U$, $F$ and $Q$ defined on all of $\R$ with compact support and maintaining the properties of the original functions. This time we carry out Step 1 for $A=(-\Delta)^s$ instead of $A=(-\Delta)^s+Q(t)$. Then by this modified Step 1 and Step 2, we get 
       \[
       \begin{split}
           0&=2\re \int_{\R}\langle  \rho_\delta\ast\eta_0 (\partial_t^2+A)U,\rho_\delta\ast \eta_0(\partial_t U)\rangle\,dt\\
               &+2\re\sum_{j=0,1}(-1)^j\langle \rho_\delta\ast \rho_\delta\ast \eta_0 \partial_t U, \eta_0\partial_t U\rangle(j t_0)\\
               &+2\re\sum_{j=0,1}(-1)^j\langle  \rho_\delta\ast \rho_\delta\ast(A\eta_0 U),\eta_0 U\rangle (j t_0),
       \end{split}
       \]
       where we used that $A$ is independent of $t$. As $U$ solves \eqref{F:eq: existence weakly viscous case}, we obtain
       \begin{equation}
        \begin{split}
            0&=-2\re\int_{\R}\langle \rho_\delta\ast  Q \eta_0 U ,\rho_\delta\ast \eta_0 (\partial_t U)\rangle\,dt\\
               &-2\re \int_{\R}\langle  \rho_\delta\ast \eta_0(-\Delta)^{r/2} \partial_t U,\rho_\delta\ast \eta_0 (-\Delta)^{r/2} \partial_t U\rangle\,dt\\
               &+2\re \int_{\R}\langle  \rho_\delta\ast F\eta_0 ,\rho_\delta\ast \eta_0(\partial_t U)\rangle\,dt\\
               &+2\re\sum_{j=0,1}(-1)^j\langle \rho_\delta\ast \rho_\delta\ast \eta_0 \partial_t U, \eta_0\partial_t U\rangle(j t_0)\\
               &+2\re\sum_{j=0,1}(-1)^j\langle  \rho_\delta\ast \rho_\delta\ast(A\eta_0 U),\eta_0 U\rangle (j t_0).
        \end{split}
       \end{equation}
       By Lemma \ref{F:lemma: cond on q} and the usual properties of mollification, we know that
       \begin{equation}
       \label{F:eq: conv delta to zero of potential term}
           \int_{\R}\langle \rho_\delta\ast  Q \eta_0 U ,(\rho_\delta\ast \eta_0 (\partial_t U))\rangle\,dt\to \int_{0}^{t_0}\langle Q  U ,\partial_t U\rangle\,dt
       \end{equation}
       as $\delta\to 0$. Thus, using the convergence results \eqref{F:eq: viscous term delta to zero}, \eqref{F:eq: source term delta to zero}, \eqref{F:eq: bdry term delta to zero}, \eqref{F:eq: remainder 1 intial cond goes to zero}, \eqref{F:eq: convergence second initial cond} and \eqref{F:eq: remainder 2 intial cond goes to zero}, we get
       \[
       \begin{split}
            0&=-2\re\int_{0}^{t_0}\langle Q  U ,\partial_t U\rangle\,dt-2 \int_{0}^{t_0}\|(-\Delta)^{r/2} \partial_t U\|_{L^2(\R^n)}^2\,dt+2\re \int_{0}^{t_0}\langle  F ,\partial_t U\rangle\,dt\\
               &+\sum_{j=0,1}(-1)^j\left(\|\partial_t U(jt_0)\|^2_{L^2(\Omega)}+\|(-\Delta)^{s/2}U(jt_0)\|^2_{L^2(\R^n)}\right).
        \end{split}
       \]
       This concludes the proof of \eqref{F:eq: energy identity weakly viscous nonlocal wave eq continuous potential}.
    \end{proof}
    \section{Proof of entanglement principle in Bessel potential spaces}
    \label{R:sec: proof entanglement principle}

\begin{proof}[Proof of Proposition~\ref{R:prop: entanglement Bessel}]
Let $\rho_\varepsilon$ be a standard mollifier and define
\[
u_{j,\varepsilon}
=
\rho_\varepsilon \ast u_j.
\]
Since $u_j\in H^{t,p}(\R^n)$, we have $u_{j,\varepsilon}\in C_b^\infty(\R^n)$ and
\[
u_{j,\varepsilon}\to u_j
\quad\text{in }H^{t,p}(\R^n)
\]
as $\varepsilon\to0$.

Moreover, convolution preserves super-exponential decay up to a modification of the constants. Indeed, arguing exactly as in the proof of
\cite[Theorem~1.2]{feizmohammadi2024entanglementprinciplefractionallaplacian},
one obtains
\[
|D^\alpha u_{j,\varepsilon}(x)|
\le
C_{\alpha,\varepsilon}
e^{-\rho |x|^\gamma}
\]
for all multiindices $\alpha$.

Since
\[
u_j=0
\quad\text{in }O,
\]
we also have
\[
u_{j,\varepsilon}=0
\quad\text{in }O_\varepsilon,
\]
where
\[
O_\varepsilon
=
\{x\in O\,;\dist(x,\partial O)>\varepsilon\}.
\]
Furthermore,
\[
\sum_{j=1}^N
b_j(-\Delta)^{s_j}u_{j,\varepsilon}
=
0
\quad\text{in }O_\varepsilon,
\]
because convolution commutes with fractional Laplacians.

Hence, the smooth entanglement principle
\cite[Theorem~1.2]{feizmohammadi2024entanglementprinciplefractionallaplacian}
implies
\[
u_{j,\varepsilon}\equiv0
\quad\text{in }\R^n
\]
for every sufficiently small $\varepsilon>0$, so that $O_\varepsilon$ contains a ball. Passing to the limit $\eps\to 0$ in $H^{t,p}(\R^n)$ yields
\[
u_j\equiv0,
\]
which proves the claim.
\end{proof}

\medskip
\noindent\textbf{Acknowledgments.} P.~Zimmermann is supported by the Swiss National Science Foundation (SNSF) under Grant Nos. 214500 and 239130. The author is particularly grateful to Katya Krupchyk for numerous stimulating and valuable discussions throughout the development of this work. Part of this research was carried out during a visit to the University of California, Irvine. The author thanks the institution for its hospitality.

\medskip 
\noindent\textbf{Declaration of AI-assisted tools.} 
Following the completion of the initial manuscript, the author used OpenAI's ChatGPT to assist with improvements to the organization and clarity of the exposition and the identification of potential proof gaps or errors. All mathematical statements and proofs were independently reviewed and verified by the author, who takes full responsibility for the content of the manuscript.

	\bibliography{refs} 
	
	\bibliographystyle{alpha}
\end{document}